\documentclass[hidelinks,onefignum,onetabnum]{siamart250106}

\definecolor{darkgreen}{rgb}{0.0, 0.55, 0.18}

\newcommand{\mycomment}[1]{%

\ifthenelse{\isodd{\value{page}}}{%
\makeatletter\normalmarginpar%
\marginpar{\tiny {#1}}%
}{%
\makeatletter\reversemarginpar%
\marginpar{\tiny {#1}}%
}}%

\usepackage[T1]{fontenc}      
\usepackage{type1cm}           
\usepackage{newtxtext}         
\usepackage[varvw]{newtxmath}  

\usepackage{graphicx}          
\usepackage{epstopdf}          
\ifpdf
  \DeclareGraphicsExtensions{.eps,.pdf,.png,.jpg}
\else
  \DeclareGraphicsExtensions{.eps}
\fi

\usepackage{amsfonts}          
\usepackage{amsopn}            

\usepackage{color}             
\usepackage{multicol}          
\usepackage{multirow}          
\usepackage{adjustbox}         
\usepackage[bottom]{footmisc}  

\usepackage{url}               
\usepackage{bm}                
\usepackage{ulem}              
\usepackage{comment}           
\usepackage{siunitx}           
\usepackage{algorithmic}       
\usepackage{makeidx}           
\usepackage{cleveref}          

\newcommand{\dd}{\mathop{}\!\mathrm{d}}

\newcommand{\Bve}{{\bm{\varepsilon}}}
\newcommand{\Bsigma}{{\bm{\sigma}}}

\newcommand{\STAB}[1]{\begin{tabular}{@{}c@{}}#1\end{tabular}}

\usepackage{acronym}           

\makeatletter
\let\AC@original@verridelabel\@verridelabel
\makeatother

\acrodef{BDDC}{Balancing Domain Decomposition by Constraints}
\acrodef{CTW}{Controlled Tensile Weldability}
\acrodef{DDMs}{domain decomposition methods}
\acrodef{DOF}{degree of freedom}
\acrodef{DOFs}{degrees of freedom}
\acrodef{FEM}{finite element methods}
\acrodef{FETI-DP}{Finite Element Tearing and Interconnecting -- Dual Primal}
\acrodef{GDSW}{Generalized Dryja-Smith-Widlund}
\acrodef{GMRES}{Generalized Minimal Residual}
\acrodef{LBW}{Laser Beam Welding}
\acrodef{PCG}{Preconditioned Conjugate Gradient}
\acrodef{PDE}{partial differential equation}
\acrodef{PDEs}{partial differential equations}
\acrodef{RGDSW}{Reduced GDSW}
\acrodef{VEM}{virtual element method}

\newsiamremark{remark}{Remark}
\newsiamremark{hypothesis}{Hypothesis}
\crefname{hypothesis}{Hypothesis}{Hypotheses}

\newsiamthm{claim}{Claim}

\headers{Two-level Schwarz Preconditioners for high-order discretizations in 3D}{T. Bevilacqua, A. Klawonn, M. Lanser and A. Wasiak}

\title{Two-level additive Schwarz preconditioners with GDSW coarse spaces for high-order discretizations in 3D}

\newcommand{\siscdeptaffil}{\thanks{Division of Mathematics, University of Cologne, Weyertal 86-90, 50931 Cologne, Germany}}
\newcommand{\sisccenteraffil}{\!\thanks{Center for Data and Simulation Science, University of Cologne, Germany}}
\newcommand{\siscemailfootnote}{\thanks{\email{E-mail: [tommaso.bevilacqua,axel.klawonn,martin.lanser,adam.wasiak]@uni-koeln.de}}}
\newcommand{\siscdeptmark}{\ensuremath{{}^{*}}}
\newcommand{\sisccentermark}{\ensuremath{\!{}^{\dagger}}}

\author{%
  Tommaso Bevilacqua\siscdeptaffil%
  \and Axel Klawonn\siscdeptmark\sisccenteraffil
  \and Martin Lanser\siscdeptmark\sisccentermark%
  \and Adam Wasiak\siscdeptmark%
  \phantom{\siscemailfootnote}.
}

\newif\ifhidesout
\newif\ifnocolors
\hidesouttrue    
\nocolorsfalse   

\ifpdf
\hypersetup{
  pdftitle={Two-level additive Schwarz preconditioners with GDSW coarse spaces for high-order discretizations in 3D},
  pdfauthor={T. Bevilacqua, A. Klawonn, M. Lanser and A. Wasiak}
}
\fi

\begin{document}

\maketitle

\begin{abstract}
Overlapping Schwarz methods are among the most successful domain decomposition preconditioners for large-scale finite element simulations. In this work, the Generalized Dryja-Smith-Widlund (GDSW), GDSW$^*$, and Reduced GDSW (RGDSW) coarse spaces are studied for three-dimensional high-order finite element methods (FEM) and virtual element methods (VEM) for the scalar diffusion equation and the equations of linear elasticity.
To the best of our knowledge, this is the first study of GDSW-type coarse spaces for polynomial degrees $k\ge 3$ in three dimensions, for either discretization.
In 3D VEM, the face degrees of freedom (DOFs) are moments instead of point values. The kernel of the differential operator, which defines the interface values of the GDSW coarse basis functions, must be represented through these moments. The construction is therefore no longer purely algebraic and, without this correction, the coarse level loses its effect for $k\ge 3$.
FEM and VEM are tested on hexahedral meshes, and VEM is additionally tested on polyhedral ones. Polynomial degrees up to $k=5$ and runs on up to $4\,096$ subdomains and parallel cores are included in the numerical study.
Numerical experiments for the scalar diffusion problem show scalable iteration counts for all coarse space variants with a constant coefficient, i.e., for the Poisson problem, while only GDSW and GDSW$^*$ remain robust for a heterogeneous high-contrast coefficient, illustrating the expected trade-off between the coarse space size and robustness.
Further experiments for linear elasticity confirm that the moment-based construction extends naturally from the scalar to the vector-valued problem.
\end{abstract}



\acresetall

\section{Introduction}
\label{intro}

Increasing the polynomial degree is a standard route to high accuracy for elliptic boundary value problems, and high-order methods are therefore an important tool in computational science and engineering~\cite{karniadakis_spectral_2005,szabo_finite_1991}.
On hexahedral grids, high-order \ac{FEM} typically relies on tensor-product spaces of type $\mathbb{Q}_k$.
The \ac{VEM} provides a different high-order framework, introduced in~\cite{beirao_basic_2013,beirao_hitchhiker_2014} and developed further in~\cite{ahmad_equivalent_2013, beirao_dassi_russo_highorder_2017, dassi_mascotto_highorder_3d_2018}, in which general polyhedral elements can be used without explicitly known basis functions.
The local \ac{VEM} functions are accessed through computable polynomial projections and through \ac{DOFs} attached to geometric entities; in 3D these include vertex and edge values as well as moments on faces.
Since its introduction, the \ac{VEM} has been applied to a wide range of problems; in particular, \ac{VEM} formulations for linear elasticity have been developed in~\cite{artioli_beirao_lovadina_sacco_vem_elasticity_2017,beirao_brezzi_marini_vem_elasticity_2013, dassi_lovadina_visinoni_vem_3d_elasticity_2020,gain_talischi_paulino_vem_elasticity_2014, xu_fan_wriggers_highorder_3d_elasticity_2024}, and recent developments include hybrid conforming/nonconforming approaches on curved three-dimensional domains~\cite{beirao_dassi_russo_trezzi_hybrid_vem_2026}.

Three-dimensional high-order FE and VE discretizations give rise to large sparse linear systems whose efficient solution requires iterative methods equipped with scalable preconditioners.
This applies to scalar problems such as the Poisson equation as well as to vector-valued problems such as linear elasticity, where the system size grows by a factor of three.
\Ac{DDMs} address this issue by coupling subdomain computations with information transferred through a global coarse problem, a mechanism that is well suited to parallel architectures~\cite{toselli_widlund_ddm_2005}.
For spectral and other high-order finite element discretizations, non-overlapping \ac{BDDC} and \ac{FETI-DP} approaches have been studied in~\cite{klawonn_pavarino_rheinbach_spectral_2008, pavarino_widlund_bddc_feti_spectral_2007, pavarino_widlund_zampini_bddc_spectral_2010}.
Within the class of non-overlapping methods for \ac{VEM}, \ac{BDDC}, and \ac{FETI-DP} preconditioners were analyzed in~\cite{bertoluzza_fetidp_2017} for two-dimensional problems and extended to three dimensions in~\cite{bertoluzza_fetidp_2020}.
Adaptive \ac{BDDC} preconditioners for high-order \ac{VEM} discretizations of three-dimensional elliptic problems in mixed form were studied in~\cite{dassi_scacchi_zampini_bddc_2022}.
More recently, \ac{BDDC} preconditioners for \ac{VEM} discretizations of the three-dimensional Stokes equations were developed in~\cite{bevilacqua_bddc_stokes_vem_2024}.
Parallel studies of high-order three-dimensional \ac{VEM} solvers have so far been limited to strong scaling on up to a few hundred processes and to mixed or saddle-point formulations~\cite{ bevilacqua_bddc_stokes_vem_2024, dassi_scacchi_parallel_mixed_2020, dassi_scacchi_zampini_bddc_2022}.
A multigrid algorithm for the $p$-version of the \ac{VEM} was proposed in~\cite{antonietti_mascotto_verani_multigrid_2018}.

For overlapping Schwarz preconditioners, classical \ac{FEM} estimates describe the dependence on the mesh and overlap parameters; see, for example,~\cite{smith_domain_1991, toselli_widlund_ddm_2005}.
For high-order discretizations, the question of robustness with respect to the polynomial degree $k$ is more delicate: in~\cite{pavarino_additive_1994}, it was shown that, for continuous high-order \ac{FEM} on structured quadrilateral meshes, a $\mathbb{Q}_1$ coarse space already yields a condition number bound independent of $k$; see also~\cite{guo_cao_preconditioner_1996, guo_cao_additive_1997} for the $hp$-version.
The authors in~\cite{antonietti_class_2011} showed that for high-order discontinuous Galerkin methods polynomial enrichment of the coarse space is needed.
When the subdomain layout is not aligned with a simple coarse mesh, \ac{GDSW} coarse spaces~\cite{dohrmann_less_regular_2008, dohrmann_family_2008} give a flexible algebraic construction.
Their basis functions are tied to interface objects, namely vertices, edges, and faces in 3D, and are obtained by discrete harmonic extension from partition-of-unity data on the interface.
The \ac{RGDSW} coarse space of~\cite{dohrmann_design_2017} keeps the same guiding principle while reducing the number of coarse unknowns, and the intermediate GDSW$^*$ variant of~\cite{hochmuth2020parallel,heinlein_monolithic_block_2026} groups vertex and edge components while retaining independent face basis functions.

Existing overlapping Schwarz analyses for \ac{VEM} have so far focused on the lowest-order case; the results of \cite{ calvo_approximation_2018,calvo_overlapping_2019} also cover high-contrast coefficients~\cite{calvo_overlapping_2024} and irregular subdomain interfaces for both \ac{VEM} and DG-VEM~\cite{calvo_virtual_2020, calvo_herrera_twolevel_discontinuous_2026}.
GDSW-type coarse spaces for three-dimensional \ac{VEM} have been considered concurrently from two directions: a coarse space for lowest-order \ac{VEM} on irregular subdomain decompositions in~\cite{aguilar_pineda_calvo_3d_irregular_2025}, and a first parallel study by the present authors covering 2D and 3D for $k=1,2$ in~\cite{bevilacqua_parallel_schwarz_2025}.
To the best of our knowledge, GDSW-type coarse spaces have not yet been studied for polynomial degrees $k\ge 3$, neither for \ac{FEM} nor for \ac{VEM}, apart from the authors'two-dimensional companion paper~\cite{wasiak_gdsw_2d_fem_vem}, which treats the Poisson problem and includes a complete condition number bound for high-order \ac{FEM} and numerical evidence of $k$-robustness for both \ac{FEM} and \ac{VEM}. For isogeometric discretizations, adaptive variants have recently been studied in~\cite{pavarino_scacchi_widlund_zampini_agdsw_iga_2026}, with high spline degrees tested only in two dimensions. In particular, the representation of the null space in the \ac{VEM} face moments, treated in the present paper, has not been addressed in the available literature.

This paper develops and tests two-level additive Schwarz preconditioners with GDSW, GDSW$^*$, and RGDSW coarse spaces for three-dimensional high-order \ac{FEM} and \ac{VEM} discretizations of two model problems: the scalar diffusion equation, which includes the Poisson equation as a special case, and the equations of linear elasticity.

The main contributions of this work are the following.
First, we show that the commonly used algebraic choice of the interface values, i.e., the vector of all ones, fails for the 3D \ac{VEM} with $k\ge 3$.
As a solution, we present a moment-based representation of the null space for the Poisson and linear elasticity problems, which applies to the GDSW, GDSW$^*$, and RGDSW coarse spaces.
This only requires the \ac{DOF} functionals of the discretization, leaving the \ac{VEM} \ac{DOFs} unchanged; the same principle applies to any discretization with non-nodal \ac{DOFs} on the interface.
With this construction, we present the first numerical study of GDSW-type coarse spaces in three dimensions for polynomial degrees up to $k=5$, for \ac{FEM} and \ac{VEM}, on hexahedral and polyhedral meshes, including heterogeneous coefficients.
Finally, we have implemented all coarse spaces in parallel within the \texttt{PETSc}-based \texttt{FE2TI} software~\cite{bevilacqua_highly_scalable_2025}, coupled with the \texttt{VEM++} library~\cite{dassi_vempp_2025}, and we report parallel scaling results on up to $4\,096$ subdomains and cores.

The remainder of the paper is structured as follows.
The two model problems and the 3D \ac{FEM} and \ac{VEM} discretizations are introduced in \mbox{\Cref{sec:discretizations}}.
The GDSW, GDSW$^*$, and RGDSW coarse spaces are presented in \mbox{\Cref{sec:GDSW}}, followed in \mbox{\Cref{sec:gdsw_high_order}} by the modifications needed for high-order \ac{VEM} face \ac{DOFs}, including the extension to the vector-valued elasticity case.
\mbox{\Cref{sec:implementation}} gives implementation details, and \mbox{\Cref{sec:results_3d}} reports the numerical results for both problems.

\section{Model problems and discretizations}
\label{sec:discretizations}

The first model problem is the scalar diffusion equation
\begin{equation}\label{eq:diffusion}
  -\nabla \cdot (\rho\,\nabla u) = f \quad \text{in } \Omega
\end{equation}
on a bounded polyhedral domain $\Omega\subset\mathbb{R}^3$, where $\rho$ is a scalar coefficient function with $0 < \rho_{\text{min}} \leq \rho \leq \rho_{\text{max}}$, for $\rho_{\text{min}},\rho_{\text{max}}\in \mathbb{R}$, taken to be constant on each element of the meshes introduced below. For $\rho\equiv 1$, \eqref{eq:diffusion} reduces to the Poisson equation $-\Delta u = f$.
We consider a partition $\partial\Omega = \Gamma_D \cup \Gamma_N$, $\Gamma_D \cap \Gamma_N = \emptyset$ with $\Gamma_D \neq \emptyset$, where Dirichlet data are prescribed on $\Gamma_D$ and a homogeneous Neumann condition $\rho\,\nabla u\cdot\mathbf{n}=0$ is imposed on $\Gamma_N$, where $\mathbf{n}$ is the outer unit normal vector to $\partial\Omega$.
Let $f\in L^2(\Omega)$ and $g\in H^{1/2}(\Gamma_D)$. Writing $u=\tilde{u}+w$, where $\tilde{u}$ is a lifting of $g$ into $H^1(\Omega)$, the unknown correction $w\in H^1_{\Gamma_D}(\Omega):=\{v\in H^1(\Omega):v|_{\Gamma_D}=0\}$ is determined by
\[
a(w,v)=F(v)\qquad \forall v\in H^1_{\Gamma_D}(\Omega),
\]
where
\[
a(w,v)=\int_\Omega \rho\, \nabla w\cdot \nabla v\,\dd x
\quad
\text{and}
\quad
F(v)=\int_\Omega f v\,\dd x-\int_\Omega \rho\, \nabla\tilde{u}\cdot\nabla v\,\dd x.
\]
We note that, since $\rho \geq \rho_{\min} > 0$, $a(v,v) = 0$ holds if and only if $\nabla v = \mathbf{0}$. Hence, in the absence of Dirichlet boundary conditions, the kernel of $a(\cdot,\cdot)$ on a connected domain consists of the constant functions, independently of $\rho$.

As a second model problem, we consider the static equations of linear elasticity for the displacement field $\mathbf{u}:\Omega\to\mathbb{R}^3$ on the same boundary decomposition.
For Lam\'e parameters $\mu>0$ and $\lambda$ with $3\lambda+2\mu>0$, we define the symmetric strain and Cauchy stress tensors
\[
\Bve(\mathbf{u})=\frac{1}{2}\bigl(\nabla\mathbf{u}+(\nabla\mathbf{u})^T\bigr)
\quad
\text{and}
\quad
\Bsigma(\mathbf{u})=2\mu\Bve(\mathbf{u})+\lambda\,\operatorname{tr}(\Bve(\mathbf{u}))\mathbf{I}.
\]
Given $\mathbf{f}\in [L^2(\Omega)]^3$ and $\mathbf{g}\in [H^{1/2}(\Gamma_D)]^3$, the weak formulation seeks $\mathbf{u}\in[H^1_{\Gamma_D}(\Omega)]^3$ such that
\[
a_{\rm el}(\mathbf{u},\mathbf{v})=F_{\rm el}(\mathbf{v})
\qquad \forall \mathbf{v}\in[H^1_{\Gamma_D}(\Omega)]^3,
\]
where
\[
a_{\rm el}(\mathbf{u},\mathbf{v})
=\int_\Omega \Bsigma(\mathbf{u}):\Bve(\mathbf{v})\,\dd x
\quad
\text{and}
\quad
F_{\rm el}(\mathbf{v})
=\int_\Omega \mathbf{f}\cdot\mathbf{v}\,\dd x.
\]

\subsection{High-order finite element discretization}
\label{sec:fem_discretization}

For the \ac{FEM} discretization we consider hexahedral meshes.
Given a tessellation $\mathcal{Q}_h$ of $\Omega$ into hexahedra with maximum element diameter $h$, the local approximation space on each $Q \in \mathcal{Q}_h$ is the tensor-product polynomial space
\begin{equation}
\mathbb{Q}_k(Q) := \operatorname{span}\{x_1^{i_1} x_2^{i_2} x_3^{i_3} : 0 \leq i_1, i_2, i_3 \leq k\}.
\end{equation}
The \ac{DOFs} are:
\begin{itemize}
    \item values at the 8 vertices of each hexahedron,
    \item for $k \geq 2$, values at $k-1$  interior Gauss--Lobatto points on each of the 12 edges,
    \item for $k \geq 2$, values at $(k-1)^2$  interior Gauss--Lobatto points on each of the 6 faces,
    \item for $k \geq 2$, values at $(k-1)^3$  interior Gauss--Lobatto points within each hexahedron.
\end{itemize}
Stiffness matrix entries are computed via Gauss--Lobatto quadrature with $k+1$ points in each direction; for the scalar diffusion problem, since $\rho$ is constant on each hexahedron $Q$, the local stiffness matrix is scaled by the value $\rho|_Q$.

\paragraph{Global space and discrete problem}
The corresponding global conforming finite element space is
\[
V_{h}^k := \bigl\{ v_h \in H^1_{\Gamma_D}(\Omega) : v_h|_{Q} \in \mathbb{Q}_k(Q) \;\; \forall Q \in \mathcal{Q}_h \bigr\}.
\]
The discrete solution of the scalar diffusion problem is the function $u_h\in V_h^k$ satisfying
\[
a(u_h,v_h)=F(v_h)\qquad \forall v_h\in V_h^k.
\]
For linear elasticity, the discrete displacement $\mathbf{u}_h\in [V_h^k]^3$ satisfies $a_{\rm el}(\mathbf{u}_h,\mathbf{v}_h)=F_{\rm el}(\mathbf{v}_h)$ for all $\mathbf{v}_h\in [V_h^k]^3$.

\subsection{Virtual element discretization in 3D}
\label{sec:vem_discretization}

The \ac{VEM} extends finite element techniques to general polyhedral meshes. We consider a polyhedral tessellation $\mathcal{T}_h$ of $\Omega$ with maximum element diameter $h$ and we follow the enhanced \ac{VEM} framework of~\cite{ahmad_equivalent_2013,beirao_dassi_russo_highorder_2017}. We denote by $\mathcal{M}_\ell(S) := \{((\mathbf{x}-\mathbf{x}_S)/h_S)^{\mathbf{s}} : |\mathbf{s}| \leq \ell\}$ the space of scaled monomials of degree at most $\ell$ on the geometric object $S$ (edge, face, or element), with centroid $\mathbf{x}_S:=\frac{1}{|S|}\int_S \mathbf{x}\,\dd S$ and diameter $h_S$, and by $\mathcal{M}^*_\ell(S)$ the subset of monomials of degree exactly $\ell$.

The construction proceeds in two steps. First, on each face $F$ of a polyhedron $K$, we construct the 2D enhanced \ac{VEM} space $V_h^k(F)$. It consists of functions in $H^1(F)$ whose traces on each edge of $F$ are polynomials of degree at most $k$, whose Laplacians belong to $\mathbb{P}_k(F)$, and which satisfy an orthogonality constraint that ensures that the $L^2$-projection $\Pi^{0,F}_k$ onto $\mathbb{P}_k(F)$, defined as 
\[
\int_F (v_h-\Pi^{0,F}_k v_h)q\,\dd x=0\qquad \forall q\in\mathbb{P}_k(F),
\]
is computable from the \ac{DOFs}. We set $V_h^k(\partial K) := \{v \in C^0(\partial K) : v|_F \in V_h^k(F) \; \forall F \subset \partial K\}$. Then, the 3D enhanced virtual element space of order $k \geq 1$ on each polyhedron $K$ is
\begin{equation}\label{eq:vem_space_3d}
\begin{aligned}
V_h^k(K) := \Big\{ v_h \in H^1(K) :\; & v_h|_{\partial K} \in V_h^k(\partial K), \quad \Delta v_h \in \mathbb{P}_k(K), \\
& (v_h - \Pi^\nabla_k v_h, \, q^*)_{0,K} = 0 \;\; \forall q^* \in \mathcal{M}^*_{k-1}(K) \cup \mathcal{M}^*_k(K) \Big\},
\end{aligned}
\end{equation}
where $\Delta$ is the standard Laplace operator, and $\Pi^\nabla_k : V_h^k(K) \to \mathbb{P}_k(K)$ is the elliptic projector defined by~\cite{beirao_dassi_russo_highorder_2017}
\begin{equation*}
\int_K \nabla(v_h - \Pi^\nabla_k v_h) \cdot \nabla p \, \dd x = 0 \quad \forall p \in \mathbb{P}_k(K),
\end{equation*}
together with the additional condition $\sum_{\nu \in \mathcal{V}(K)} (v_h(\nu) - \Pi^\nabla_k v_h(\nu)) = 0$ for $k = 1$, and $\int_K (v_h - \Pi^\nabla_k v_h) \, \dd x = 0$ for $k \geq 2$.
This space satisfies $\mathbb{P}_k(K) \subseteq V_h^k(K)$.
Following~\cite{beirao_dassi_russo_highorder_2017}, the \ac{DOFs} for $V_h^k(K)$ are:
\begin{itemize}
    \item values of $v_h$ at the vertices of $K$,
    \item for $k \geq 2$, values of $v_h$ at the $k-1$ internal Gauss--Lobatto points on each edge of $\partial K$,
    \item for $k \geq 2$, face moments on each face $F$ of $K$:
    $\frac{1}{|F|}\int_F p \, v_h \, \dd S \quad \text{for all } p \in \mathcal{M}_{k-2}(F),$
    \item for $k \geq 2$, interior moments in $K$:
    $\frac{1}{|K|}\int_K q \, v_h \, \dd x \quad \text{for all } q \in \mathcal{M}_{k-2}(K)$.
\end{itemize}
The use of the 2D enhanced space $V_h^k(F)$ on each face is essential in 3D. It allows for the computation of the full $L^2$-projection $\Pi^{0,F}_k$ onto $\mathbb{P}_k(F)$, which is needed to compute the boundary integral in the construction of the 3D projector $\Pi^\nabla_k$~\cite{ahmad_equivalent_2013}.
The projector $\Pi^\nabla_k$ enters the discrete form through a consistency term on the projected polynomials and a stabilization on the remainder,
\[
a_h(v_h,w_h)=a(\Pi^\nabla_k v_h,\Pi^\nabla_k w_h)+s_h(v_h-\Pi^\nabla_k v_h,w_h-\Pi^\nabla_k w_h).
\]
Since $\rho$ is constant on each element $K$, the projector $\Pi^\nabla_k$ does not depend on $\rho$, and on $K$, both, the consistency term and the stabilization term $s_h$, are scaled by the value $\rho|_K$.
The global conforming virtual element space is
\[
V_{h}^k := \bigl\{ v_h \in H^1_{\Gamma_D}(\Omega) : v_h|_{K} \in V_h^k(K) \;\; \forall K \in \mathcal{T}_h \bigr\}.
\]
We introduce the computable $L^2$-projector $\Pi^0_k$ obtained from the defining relation
\[
\int_K (v_h-\Pi^0_k v_h)q\,\dd x=0\qquad \forall q\in\mathbb{P}_k(K).
\]
The discrete right-hand side is defined through $F_h(v_h):=F(\Pi^0_k v_h)$, and the resulting \ac{VEM} solution is the function $u_h\in V_h^k$ satisfying $a_h(u_h,v_h)=F_h(v_h)$ for every $v_h\in V_h^k$.
For the standard assumptions and convergence results, see~\cite{ahmad_equivalent_2013,beirao_basic_2013,beirao_stability_2017}.

\paragraph{VEM discretization for linear elasticity}
The local virtual element space for linear elasticity is the vector-valued space $[V_h^k(K)]^3$.
The \ac{DOFs} are the three displacement components associated with each scalar DOF: vertex values, edge Gauss--Lobatto values, face moments, and element interior moments.
The stress and the strain are approximated by their $L^2$-projections onto polynomials of degree $k-1$, denoted by $\boldsymbol{\Pi}^0_{k-1}$, while the elliptic projector $\boldsymbol{\Pi}_k^\nabla$ is applied componentwise in the stabilization; on each element $K$, the discrete bilinear form is
\[
a_{h,\rm el}(\mathbf{v}_h,\mathbf{w}_h)
= \int_K \boldsymbol{\Pi}^0_{k-1}\Bsigma(\mathbf{v}_h):\boldsymbol{\Pi}^0_{k-1}\Bve(\mathbf{w}_h)\,\dd x
+ \mathbf{s}_h(\mathbf{v}_h - \boldsymbol{\Pi}_k^\nabla \mathbf{v}_h,\, \mathbf{w}_h - \boldsymbol{\Pi}_k^\nabla \mathbf{w}_h),
\]
where $\mathbf{s}_h(\cdot,\cdot)$ is the componentwise d-recipe stabilization.

\section{Overlapping Schwarz method with GDSW}
\label{sec:GDSW}

Both the \ac{FEM} and \ac{VEM} discretizations of \mbox{\Cref{sec:discretizations}} yield a sparse linear system
\begin{equation}\label{eq:linear_system}
A\mathbf{u} = \mathbf{f},
\end{equation}
where $A\in\mathbb{R}^{n\times n}$ is the assembled stiffness matrix, $\mathbf{u}\in\mathbb{R}^n$ is the vector of \ac{DOF} values, and $\mathbf{f}\in\mathbb{R}^n$ is the right-hand side.
For the scalar diffusion problem, $\mathbf{u}$ collects one component per \ac{DOF}; for linear elasticity, it collects three displacement components, so $n$ is three times the number of scalar \ac{DOFs}.
The algebraic structure of the preconditioner and the coarse space construction described below apply to both problems: only the definition of the DOF vector and the operator $A$ changes, while the interior--interface splitting, the discrete harmonic extension, and the partition-of-unity interface data are defined identically.
For ease of presentation, the GDSW construction is described for the scalar diffusion problem; the extension to the vector-valued elasticity case, which requires specific modifications, is addressed in \mbox{\Cref{sec:gdsw_high_order}}.
We use a two-level overlapping Schwarz preconditioner for~\eqref{eq:linear_system} whose global coarse correction is spanned by GDSW-type coarse basis functions~\cite{dohrmann_less_regular_2008,dohrmann_family_2008}. The local solvers are defined on the subdomains of an overlapping decomposition $\Omega = \bigcup_{i=1}^N \Omega_i$, where the subdomains are usually obtained by adding $\delta$ rows of finite or virtual elements to non-overlapping ones.
Now, the two-level additive Schwarz operator consists of the sum of the coarse correction and local solvers,
\[
M^{-1}
=\Phi A_0^{-1}\Phi^T+\sum_{i=1}^N R_i^T A_i^{-1}R_i^{\vphantom{-1}}.
\]
Here, $R_i$ restricts to the overlapping subdomain $\Omega_i$, $A^{\vphantom{T}}_i=R^{\vphantom{T}}_i AR_i^T$, and the columns of $\Phi$ span the chosen coarse space, with the coarse matrix $A_0=\Phi^T A\Phi$.
With a slight abuse of notation, we use $M^{-1}$ for all coarse space variants introduced later; the specific choice of coarse basis functions collected in $\Phi$ will be clear from the context of the numerical experiments.
For the three-dimensional GDSW space, the non-overlapping interface is split into vertices, edges, and faces.
The interface components are identified algebraically. For each \ac{DOF}, we determine the set of non-overlapping subdomains that contains it; the \ac{DOF} is an interface \ac{DOF} if it belongs to at least two subdomains. Interface \ac{DOFs} with the same set of subdomains form an equivalence class, and each class is one interface component. Specifically, a class of a single \ac{DOF} is a vertex, while a class of several \ac{DOFs} is a face if it is shared by exactly two subdomains and an edge otherwise.
We denote by $\mathcal{P}$ a fixed interface component, which we write as $\mathcal{V}$, $\mathcal{E}$, or $\mathcal{F}$, respectively, when the specific type needs to be indicated. We use $\Gamma$ for the entire non-overlapping interface.
Each such interface component $\mathcal{P}$ defines boundary data $\phi_{\mathcal{P}}$, a vector defined on the whole interface $\Gamma$ that vanishes outside $\mathcal{P}$. We then accordingly write $\phi_{\mathcal{V}}$, $\phi_{\mathcal{E}}$, and $\phi_{\mathcal{F}}$ for the boundary data of a vertex, edge, and face component, respectively.
This partition introduces a splitting of the \ac{DOFs} in the global stiffness matrix $A$ into interior ($I$) and interface ($\Gamma$),
\[
A = \begin{pmatrix} A_{II} & A_{I\Gamma} \\ A_{\Gamma I} & A_{\Gamma\Gamma} \end{pmatrix},
\]
where $A_{II}$ is the submatrix corresponding to interior \ac{DOFs} and $A_{I\Gamma}$ couples interior to interface \ac{DOFs}. The interior values of the coarse basis functions are then obtained by the discrete harmonic extension
\begin{equation}\label{eq:harmonic_ext}
\phi_I^{\mathcal{P}}=- A_{II}^{-1} A^{\vphantom{-1}}_{I\Gamma}\phi^{\vphantom{-1}}_{\mathcal{P}}.
\end{equation}
In~\cite{dohrmann_family_2008}, the interface data are obtained by restricting the null space of the operator to the interface components. For scalar problems, the null space consists of the constant functions, so that $\phi_{\mathcal{P}}$ contains the DOF values of the constant function $1$ on $\mathcal{P}$ and vanishes on the rest of the interface. For nodal \ac{DOFs}, as in the scalar examples of~\cite{dohrmann_family_2008}, these values are all equal to one.
Vertex and edge \ac{DOFs} are always nodal: $\phi_{\mathcal{V}}(x)=1$ for $x\in\mathcal{V}$ and $0$ otherwise, and analogously for $\phi_{\mathcal{E}}$. For \ac{FEM}, $\phi_{\mathcal{F}}=1$ on all face \ac{DOFs} belonging to $\mathcal{F}$ and $0$ elsewhere, while for high-order \ac{VEM} in three dimensions the assignment requires a moment-based evaluation described in \mbox{\Cref{sec:gdsw_high_order}}. For each vertex, edge, and face, the corresponding coarse basis function is then defined as $\phi^{\mathcal{P}}=((\phi_I^{\mathcal{P}})^T,\phi^{T}_{\mathcal{P}})^T$, which is one column of $\Phi$. Denoting by $n_\mathcal{V}$, $n_\mathcal{E}$, and $n_\mathcal{F}$ the number of vertices, edges, and faces of the non-overlapping interface decomposition, the standard \ac{GDSW} coarse space is thus spanned by $\Phi\in\mathbb{R}^{n\times N_c}$ with $N_c=n_\mathcal{V}+n_\mathcal{E}+n_\mathcal{F}$.

\subsection{RGDSW}
\label{sec:rgdsw}

The \ac{RGDSW} coarse space~\cite{dohrmann_design_2017} is a variant of the standard \ac{GDSW} coarse space designed to reduce the dimension of the coarse problem while maintaining robust convergence properties. In this approach, the separate edge and face coarse functions are replaced by vertex-centered combinations.
For a vertex $\mathbf{v}$, let $\mathcal{E}_{\mathbf{v}}$ and $\mathcal{F}_{\mathbf{v}}$ collect the incident interface edge and face \ac{DOFs}, and let $n_{\mathcal{E}_{\mathbf{v}}}$ and $n_{\mathcal{F}_{\mathbf{v}}}$ be the number of non-Dirichlet vertices adjacent to the edge or face, respectively.
The interface support of the RGDSW basis function and the interface values are given by
\begin{equation}\label{eq:rgdsw_basis}
\mathcal{P}_{\mathbf{v}}=\{\mathbf{v}\}\cup\mathcal{E}_{\mathbf{v}}\cup\mathcal{F}_{\mathbf{v}}, \qquad \text{and} \qquad \phi_{\mathcal{P}_{\mathbf{v}}}(x)=
\begin{cases}
1 & x=\mathbf{v}, \\
\tfrac{1}{n_{\mathcal{E}_\mathbf{v}}} & x\in \mathcal{E}_{\mathbf{v}}, \\
\tfrac{\phi_{\mathcal{F}_\mathbf{v}}}{n_{\mathcal{F}_\mathbf{v}}} & x\in \mathcal{F}_{\mathbf{v}}, \\
0 & \text{else},
\end{cases}
\end{equation}
respectively, where $\phi_{\mathcal{F}_\mathbf{v}}$ is the face value defined as in \mbox{\Cref{sec:GDSW}}.
The remaining entries are computed by discrete harmonic extension~\eqref{eq:harmonic_ext}, giving the full coarse vector $\phi^{\mathcal{P}_{\mathbf{v}}}=((\phi_I^{\mathcal{P}_{\mathbf{v}}})^T,\phi^T_{\mathcal{P}_{\mathbf{v}}})^T$.

\subsection{GDSW$^*$}
\label{sec:gdswstar}

The GDSW$^*$ coarse space, introduced in~\cite{hochmuth2020parallel,heinlein_monolithic_block_2026} and further discussed in~\cite{bevilacqua_parallel_schwarz_2025}, is an intermediate variant between GDSW and RGDSW.
While RGDSW groups all interface components (vertices, edges, and faces) into a single basis function per vertex, GDSW$^*$ groups only the vertex and edge components, retaining one independent basis function per face as in the standard GDSW construction.
This yields $n_\mathcal{V} + n_\mathcal{F}$ coarse basis functions, compared to $n_\mathcal{V} + n_\mathcal{E} + n_\mathcal{F}$ for GDSW and $n_\mathcal{V}$ for RGDSW.

The vertex basis function $\phi^{\mathcal{P}_{\mathbf{v}},*}$ is defined by prescribing
\begin{equation}
\label{eq:gdswstar_basis}
\phi^{*}_{\mathcal{P}_{\mathbf{v}}}(x) = \begin{cases}
1 & x = \mathbf{v}, \\
\tfrac{1}{n_{\mathcal{E}_\mathbf{v}}} & x \in \mathcal{E}_{\mathbf{v}}, \\
0 & \text{else},
\end{cases}
\end{equation}
and the interior part is obtained by discrete harmonic extension~\eqref{eq:harmonic_ext}.
The face basis functions are identical to those of the standard GDSW construction, with $\phi_{\mathcal{F}}$ as defined in \mbox{\Cref{sec:GDSW}}.

\section{GDSW-type coarse spaces for high-order \ac{VEM}}
\label{sec:gdsw_high_order}

The decisive step in constructing GDSW-type coarse spaces for high-order \ac{VEM} is the choice of interface data before the harmonic extension is computed.
For the scalar diffusion operator, the coarse correction must at least reproduce constants, since they represent its kernel for any $\rho$ that satisfies our assumptions (cf.\ \mbox{\Cref{sec:discretizations}}).
Thus, the 3D construction has to assign vertex, edge, and face \ac{DOFs} such that the assembled coarse functions still form a partition of unity on the interface, i.e., that the prescribed boundary data correctly represent the constant function in the \ac{DOFs} of the discrete space.

For nodal \ac{DOFs}, i.e., for all \ac{FEM} \ac{DOFs} and for the vertex and edge \ac{DOFs} of the \ac{VEM} of the discretizations used in this work, this is achieved by the values $1$ and $0$ of \mbox{\Cref{sec:GDSW}}. The faces of the subdomain interface, however, also carry the \ac{VEM} face moments
\begin{equation}
\frac{1}{|F|}\int_F p \, v_h \, \dd S \quad \text{for all } p \in \mathcal{M}_{k-2}(F),
\end{equation}
cf. \mbox{\Cref{sec:vem_discretization}}. Assigning the value 1 to all such moments would incorrectly represent the constant function on the interface, since the moments of the constant function tested against high-order monomials are generally not equal to 1.

Therefore, for each subdomain face $\mathcal{F}$, we define the face DOF values of the GDSW coarse basis function $\phi_\mathcal{F}$ by evaluating the face moments against the constant function $1 \in \mathbb{P}_k(F) \subseteq V_h^k(F)$. Specifically, for each element face $F \subset \mathcal{F}$, let $m_1, m_2, \ldots$ be an enumeration of the basis monomials of $\mathcal{M}_{k-2}(F)$. We denote by $(\phi_\mathcal{F})_{F,j}$ the value of the $j$-th face DOF of $\phi_\mathcal{F}$ on element face $F$, and set
\begin{equation}\label{eq:gdsw_vem_face}
(\phi_\mathcal{F})_{F,j} := \frac{1}{|F|}\int_{F} m_j \, \dd S, \quad j = 1, \ldots, \dim(\mathcal{M}_{k-2}(F)), \quad \forall\, F \subset \mathcal{F}.
\end{equation}
The nodal \ac{DOFs} of the mesh vertices and mesh edges lying in $\mathcal{F}$ are set to $1$ as before, and all other interface \ac{DOFs} are set to $0$.
The RGDSW and GDSW$^*$ coarse spaces follow the constructions of \mbox{\Cref{sec:rgdsw,sec:gdswstar}}, respectively, with the face moment \ac{DOFs} computed via~\eqref{eq:gdsw_vem_face} in place of the standard nodal assignment. This is the choice required by the definition in \mbox{\cite{dohrmann_family_2008}} as we recall in \mbox{\Cref{sec:GDSW}}. There, however, non-nodal \ac{DOFs} appear only for lowest-order edge elements, with the vector of all ones and a suitable scaling of the basis functions; high-order moments are not considered.
\begin{remark}
The moment-based assignment~\eqref{eq:gdsw_vem_face} differs from the standard nodal assignment for $k\ge 3$.
For $k=1$, no face moments are present and all interface \ac{DOFs} are nodal, so the usual and modified constructions coincide.
For $k=2$, face moments appear but the test space is $\mathcal{M}_{0}(F)$, which consists of the constant monomial alone; the corresponding moment value is $\frac{1}{|F|}\int_F 1\,\dd S = 1$, which is exactly the value assigned by the usual construction.
Hence the two constructions produce identical coarse spaces for $k\le 2$.
For $k=3$, the test space is $\mathcal{M}_{1}(F)=\big\{1,\,(\mathbf{x}-\mathbf{x}_F)_1/h_F,\,(\mathbf{x}-\mathbf{x}_F)_2/h_F\big\}$, where the components refer to the local coordinates on $F$. The moment of the constant monomial is again $1$, while
\[
\frac{1}{|F|}\int_F \frac{(\mathbf{x}-\mathbf{x}_F)_i}{h_F}\,\dd S = \frac{1}{h_F}\Big(\frac{1}{|F|}\int_F \mathbf{x}_i\,\dd S-(\mathbf{x}_F)_i\Big) = 0,\qquad i=1,2,
\]
by the definition of the centroid. The face moment \ac{DOFs} of the constant function on each element face are thus given by $(1,0,0)$, whereas the standard assignment would instead result in $(1,1,1)$.
For $k\ge 4$, the moments of the monomials of degree two and higher are in general neither $0$ nor $1$ and depend on the shape of $F$.
\end{remark}
\begin{remark}
Face moments appear on the subdomain interface because faces are part of the non-overlapping decomposition boundary in 3D. The same consideration applies to any discretization method involving non-nodal \ac{DOFs} on interface components; the DOF values for the coarse basis must be set to the values of the \ac{DOF} functionals applied to the functions spanning the kernel of the differential operator, rather than to 1. Since these values have to be provided by the discretization, the construction of the coarse space is no longer purely algebraic.
\end{remark}

\paragraph{Extension to linear elasticity}
For the linear elasticity operator $-\nabla\cdot\Bsigma(\cdot)$ (cf.\ \mbox{\Cref{sec:discretizations}}), the kernel consists of the six rigid-body modes: the three constant translations
\[
\mathbf{r}_1=\begin{pmatrix}1\\0\\0\end{pmatrix},\quad
\mathbf{r}_2=\begin{pmatrix}0\\1\\0\end{pmatrix},\quad
\mathbf{r}_3=\begin{pmatrix}0\\0\\1\end{pmatrix},
\]
and the three rotations
\[
\mathbf{r}_4=\begin{pmatrix}0\\x_3\\-x_2\end{pmatrix},\quad
\mathbf{r}_5=\begin{pmatrix}-x_3\\0\\x_1\end{pmatrix},\quad
\mathbf{r}_6=\begin{pmatrix}x_2\\-x_1\\0\end{pmatrix}.
\]
All six satisfy $\Bve(\mathbf{r}_i)=0$, hence $a_{\rm el}(\mathbf{r}_i,\mathbf{v})=0$ for all $\mathbf{v}$ and $i=1,\ldots,6$.
The GDSW coarse space for elasticity contains up to one coarse basis function per interface component and rigid-body mode, so that the assembled coarse functions reproduce all six rigid-body modes on the interface.
As usual, rigid-body modes that are linearly dependent on an interface component are discarded: a vertex carries only the three translations, since at a single node the rotations are linear combinations of them, and a straight edge carries five coarse basis functions, since the rotation about the edge direction reduces to a translation on the edge. All other interface components have six rigid-body modes.

Since the vector-valued \ac{VEM} space $[V_h^k(K)]^3$ uses the same moment-based face \ac{DOFs} as the scalar space, one copy per displacement component, the face assignment generalizes~\eqref{eq:gdsw_vem_face} by testing against each rigid-body mode.
For each face $\mathcal{F}$, each rigid-body mode $\mathbf{r}_i$ ($i=1,\ldots,6$), each element face $F\subset\mathcal{F}$, and each displacement component $\beta=1,2,3$, the $j$-th face moment DOF of the coarse basis function $\phi_\mathcal{F}^i$ is
\begin{equation}\label{eq:gdsw_vem_face_el}
(\phi_\mathcal{F}^i)_{F,j,\beta} := \frac{1}{|F|}\int_{F} m_j \, (\mathbf{r}_i)_\beta \, \dd S, \quad j = 1, \ldots, \dim(\mathcal{M}_{k-2}(F)).
\end{equation}
For the translations $\mathbf{r}_1,\mathbf{r}_2,\mathbf{r}_3$, one has $(\mathbf{r}_i)_\beta=\delta_{i\beta}$, so~\eqref{eq:gdsw_vem_face_el} reduces to the scalar formula~\eqref{eq:gdsw_vem_face} on the active component and vanishes on the others.
For the rotations $\mathbf{r}_4,\mathbf{r}_5,\mathbf{r}_6$, the components $(\mathbf{r}_i)_\beta$ are linear polynomials in $\mathbf{x}$, so the product $m_j(\mathbf{r}_i)_\beta$ is a polynomial and the integral in~\eqref{eq:gdsw_vem_face_el} is computable from the \ac{DOFs}.
Vertex and edge \ac{DOFs} on $\mathcal{F}$ are set to $(\mathbf{r}_i)_\beta(\mathbf{x})$ at each DOF coordinate value $\mathbf{x}$; all other interface \ac{DOFs} are set to zero.
The RGDSW and GDSW$^*$ coarse spaces for elasticity follow the same structural constructions as in the scalar case (see \mbox{\Cref{sec:rgdsw,sec:gdswstar}}) with the interface values described above.
As in the scalar case, the moment-based assignment coincides with the standard nodal one for $k\le 2$: for $k=1$ only vertex \ac{DOFs} are present, while for $k=2$ the face moments test against the constant monomial and the integral~\eqref{eq:gdsw_vem_face_el} reduces to the mean value of $(\mathbf{r}_i)_\beta$ on $F$, which for the linear rigid-body modes coincides with the evaluation at the face centroid.
For $k\ge 3$, the test space $\mathcal{M}_{k-2}(F)$ includes higher-order monomials whose moments against the rigid-body modes are generally non-trivial, and the moment-based assignment~\eqref{eq:gdsw_vem_face_el} becomes necessary.

\subsection{Remarks on theoretical analysis}
\label{sec:theory_remarks}
In the companion paper \cite{wasiak_gdsw_2d_fem_vem}, the authors provide a theoretical analysis of the GDSW and RGDSW preconditioners for high-order \ac{FEM} in two dimensions. The analysis establishes a condition number bound of the form
\begin{equation*}
        \kappa(M^{-1}A) \leq C \Bigl(1 + \tfrac{H}{\delta}\Bigr)\Bigl(1+\log\Bigl(\tfrac{H}{h}\Bigr)\Bigr)\Bigl(1+\log\Bigl(\tfrac{Hk^2}{h}\Bigr)\Bigr),
    \end{equation*}
for $M^{-1}\in\{M_\mathrm{GDSW}^{-1}, M_\mathrm{RGDSW}^{-1}\}$.
Currently, to the authors' knowledge, the development of technical tools, such as a face lemma or a discrete Sobolev inequality, required for a similar analysis of GDSW-type preconditioners in three dimensions, has been restricted to the $k=1$ case \cite{dryja_smith_widlund_schwarz_1994,klawonn_widlund_dryja_dual_primal_2002,toselli_widlund_ddm_2005}.
An exception is the case of spectral elements, for which such tools have been developed for iterative substructuring methods in~\cite{casarin_quasi_optimal_1997,pavarino_widlund_polylogarithmic_1996}, and have not yet been used for GDSW-type coarse spaces.
Some noteworthy results on high-order additive Schwarz methods in three dimensions can be found in \cite{guo_cao_additive_1998}.

\section{Implementation}
\label{sec:implementation}

All Schwarz preconditioners and GDSW-type coarse spaces described in this paper are incorporated into the \texttt{FE2TI} software package~\cite{bevilacqua_highly_scalable_2025,bevilacqua_klawonn_lanser_wasiak_2024}, a C/C++ library employing a hybrid MPI/OpenMP parallelization strategy built on top of the \texttt{PETSc} library~\cite{balay_petsc_2024}. In addition to the overlapping Schwarz preconditioners used here, \texttt{FE2TI} includes efficient implementations of the non-overlapping domain decomposition methods \ac{BDDC} and \ac{FETI-DP}, for which it has demonstrated excellent strong and weak scalability, handling simulations with millions of parallel ranks efficiently~\cite{klawonn_koehler_lanser_rheinbach_2020,klawonn_lanser_rheinbach_uran_2021}. A detailed description of the software infrastructure regarding the Schwarz preconditioners can be found in~\cite{bevilacqua_highly_scalable_2025}.

Within the preconditioner, the first level uses a symmetric additive Schwarz method with an overlap of one layer of  elements. The symmetric formulation is required for the PCG solver and allows for the estimation of condition numbers using the Lanczos process. A restricted additive Schwarz variant~\cite{cai_restricted_1999} can reduce communication overhead and improve convergence speed in practice, but it breaks the symmetry of the preconditioned operator; one would then need to replace PCG with a non-symmetric Krylov method such as GMRES. This is available in \texttt{FE2TI}, but since a goal of this work is the analysis of condition numbers, we use the symmetric additive formulation and PCG. The second level uses a GDSW-type coarse space, with the specific variant (GDSW, RGDSW, or GDSW$^*$) chosen according to the numerical experiment.

All local subdomain problems are solved by PARDISO from the Intel Math Kernel Library (MKL)~\cite{intel_onemkl_pardiso,schenk2001pardiso}. The coarse basis functions are constructed by solving local Dirichlet extension problems and assembled into an MPI-parallel prolongation matrix $\Phi$, and the coarse stiffness matrix $A_0$ is obtained via the \texttt{PtAP} product from \texttt{PETSc}. To ensure an efficient factorization, $A_0$ is redistributed to a reduced number of MPI ranks using a dedicated sub-communicator and then factored by MKL Cluster PARDISO.

For the \ac{VEM} experiments, the element-level stiffness matrices, right-hand sides, and the moment values needed for the coarse space construction of \mbox{\Cref{sec:gdsw_high_order}} are computed by \texttt{VEM++}~\cite{dassi_vempp_2025} and passed to \texttt{FE2TI} via a dedicated interface. All subsequent steps, i.e., global parallel assembly, coarse space construction, and iterative solve, are handled entirely by \texttt{FE2TI}. This coupling allows \texttt{VEM++} to serve as a black-box local assembler and to exploit \texttt{FE2TI}'s optimized parallel infrastructure. The polynomial basis in the \ac{VEM} projectors follows the scaled monomial construction of~\cite{beirao_dassi_russo_highorder_2017}, as implemented in \texttt{VEM++}.
For the stabilization term in the \ac{VEM} bilinear form, we use the d-recipe stabilization of~\cite{beirao_dassi_russo_highorder_2017} with scaling weight equal to $1$ for both problems: for the scalar diffusion problem the stabilization uses the diagonal entries of the consistency part of the \ac{VEM} stiffness, while for the linear elasticity problem the same d-recipe is applied componentwise to each displacement component (cf.\ \mbox{\Cref{sec:discretizations}}).
In the \texttt{VEM++} library, this corresponds to the \texttt{RECIPESTAB} stabilization.

\section{Numerical results}
\label{sec:results_3d}

The computational study considers two geometric configurations and both model problems introduced in \mbox{\Cref{sec:discretizations}}.
For the scalar diffusion problem \eqref{eq:diffusion}, we set $\rho\equiv 1$, i.e., we solve the Poisson equation, in all experiments except the four-channel problem in \cref{sec:results_channel}, where $\rho$ is heterogeneous.
The first three experiments in \mbox{\Cref{sec:results_3d_dd,sec:results_3d_strong,sec:results_3d_scal}} use the unit cube $\Omega=(0,1)^3$ with homogeneous Dirichlet boundary conditions on the entire boundary and the right-hand side defined in order to have the exact solution
\[
u(x,y,z)=\sin(\pi x)\sin(\pi y)\sin(\pi z).
\]
The experiments in \cref{sec:results_3d_strip,sec:results_3d_tol,sec:results_channel} use an elongated slab domain $\Omega=(0,64)\times(0,1)\times(0,1)$ with homogeneous Dirichlet boundary conditions on the faces $\{x=0\}$ and $\{x=64\}$, homogeneous Neumann boundary conditions on the four remaining faces, and the constant right-hand side $f=1$.
The final experiment in \cref{sec:results_elasticity} validates the moment-based construction for the linear elasticity problem on the slab $\Omega=(0,32)\times(0,1)\times(0,1)$, with the face $\{x=0\}$ clamped and the remaining faces free.
All linear systems are solved by PCG, where the convergence criterion is a reduction of $10^{-8}$ of the relative natural norm  $\|M^{-1/2}r\|/\|M^{-1/2}r_0\|$; here, $\| \cdot \|$ denotes the $\ell_2$-norm. This corresponds to \texttt{PETSc}'s \texttt{-ksp\_norm\_type natural} setting, which tracks the energy norm of the error and is the naturally available norm within the CG iteration without requiring any extra inner products.
The reported condition numbers are Lanczos estimates also available in \texttt{PETSc} through the \texttt{-ksp\_monitor\_singular\_value} option.
The maximum number of PCG iterations is set to $1000$; when this limit is reached, the corresponding entry is marked N.C.\ (not converged) in the tables.
In all experiments the overlap is set to one layer of elements ($\delta=1$).

Throughout the tables in this section, ASM~(1-lv) denotes the one-level additive Schwarz method, while GDSW, GDSW$^*$, and RGDSW refer to the two-level methods employing the respective coarse space variants constructed as described in \mbox{\Cref{sec:GDSW,sec:gdswstar}}, with the moment-based construction for \ac{VEM}. The column labeled "Coarse" reports the dimension of the coarse space.
For the structured hexahedral meshes, this dimension is independent of $k$. 
For the CVT meshes partitioned by METIS, instead, the subdomain boundaries are irregular and some interface components consist of very few \ac{DOFs}, so that the classification of \mbox{\Cref{sec:GDSW}} can change with $k$. A typical case is two vertices joined by a single mesh edge, whose sets of subdomains share more than two subdomains. For $k=1$, the mesh edge has no \ac{DOF}, and, if no other \ac{DOF} is shared by exactly the common subdomains, there is no interface edge between the two vertices. For $k\ge 2$, the $k-1$ Gauss--Lobatto \ac{DOFs} on the mesh edge are shared exactly by the common subdomains and form an additional interface component: an extra edge for $k\ge 3$, and an extra vertex for $k=2$, where it consists of a single \ac{DOF}. The same happens for faces through the \ac{DOFs} on mesh edges and mesh faces, so that the coarse space dimension slightly increases with $k$; since there is not a substantial difference, we always report the value for $k=1$.

\subsection{Domain decomposition strong scaling}
\label{sec:results_3d_dd}

\begin{table}[tbp]
\centering
\caption{Condition numbers and PCG iteration counts for $\mathcal{Q}_k$-FEM and \ac{VEM} on a $48\times 48\times 48$ hexahedral mesh. \ac{DOFs} for \ac{FEM}: $117\,649$, $912\,673$, $3\,048\,625$, $7\,189\,057$, $13\,997\,521$ for $k=1,\ldots,5$; \ac{DOFs} for \ac{VEM}: $117\,649$, $912\,673$, $2\,267\,569$, $4\,292\,929$, $7\,099\,345$ for $k=1,\ldots,5$.}\label{tab:3d_160_kappa}
\begin{adjustbox}{width=\textwidth}
\begin{tabular}{|c|c|r|r||r|r||r|r||r|r||r|r||r|r|}
\cline{3-14}
\multicolumn{2}{c|}{} & nSub & Coarse & \multicolumn{2}{c||}{$k=1$} & \multicolumn{2}{c||}{$k=2$} & \multicolumn{2}{c||}{$k=3$} & \multicolumn{2}{c||}{$k=4$} & \multicolumn{2}{c|}{$k=5$} \\
\cline{3-14}
\multicolumn{2}{c|}{} &  &  & $\kappa$ & it & $\kappa$ & it & $\kappa$ & it & $\kappa$ & it & $\kappa$ & it \\
\cline{3-14}
\noalign{\vskip 0.6ex}
\cline{1-14}
\multirow{12}{*}{\STAB{\rotatebox[origin=c]{90}{FEM}}}
& \multirow{4}{*}{\STAB{\rotatebox[origin=c]{90}{GDSW}}}
& $4\times 4\times 4$ & 279   & { 26.09} & { 32} & { 39.39} & { 35} & { 47.21} & { 39} & { 52.52} & { 40} & { 56.35} & { 40} \\
& & $6\times 6\times 6$ & 1115   & { 21.34} & { 29} & { 34.05} & { 36} & { 41.84} & { 39} & { 47.41} & { 40} & { 51.62} & { 41} \\
& & $8\times 8\times 8$ & 2863   & { 17.99} & { 28} & { 30.22} & { 35} & { 37.99} & { 38} & { 43.70} & { 39} & { 48.10} & { 40} \\
& & $12\times 12\times 12$ & 10439   & { 13.42} & { 26} & { 24.70} & { 33} & { 32.36} & { 37} & { 38.15} & { 39} & { 42.70} & { 40} \\
\cline{2-14}
& \multirow{4}{*}{\STAB{\rotatebox[origin=c]{90}{GDSW$^*$}}}
& $4\times 4\times 4$ & 171   & { 25.21} & { 32} & { 37.46} & { 36} & { 44.48} & { 39} & { 49.25} & { 39} & { 52.69} & { 40} \\
& & $6\times 6\times 6$ & 665   & { 21.73} & { 31} & { 33.63} & { 36} & { 40.56} & { 39} & { 45.53} & { 40} & { 49.31} & { 42} \\
& & $8\times 8\times 8$ & 1687   & { 18.91} & { 30} & { 30.67} & { 36} & { 37.61} & { 39} & { 42.73} & { 41} & { 46.67} & { 42} \\
& & $12\times 12\times 12$ & 6083   & { 14.11} & { 27} & { 25.82} & { 34} & { 32.99} & { 38} & { 38.33} & { 41} & { 42.51} & { 42} \\
\cline{2-14}
& \multirow{4}{*}{\STAB{\rotatebox[origin=c]{90}{RGDSW}}}
& $4\times 4\times 4$ & 27   & { 33.85} & { 30} & { 50.57} & { 34} & { 60.24} & { 37} & { 66.53} & { 39} & { 70.91} & { 40} \\
& & $6\times 6\times 6$ & 125   & { 28.51} & { 32} & { 44.39} & { 37} & { 53.83} & { 40} & { 60.38} & { 42} & { 65.19} & { 43} \\
& & $8\times 8\times 8$ & 343   & { 24.14} & { 31} & { 39.40} & { 39} & { 48.56} & { 42} & { 55.08} & { 44} & { 60.00} & { 45} \\
& & $12\times 12\times 12$ & 1331   & { 17.06} & { 29} & { 31.71} & { 37} & { 40.71} & { 41} & { 47.21} & { 44} & { 52.17} & { 45} \\
\hline
\hline
\multirow{12}{*}{\STAB{\rotatebox[origin=c]{90}{VEM}}}
& \multirow{4}{*}{\STAB{\rotatebox[origin=c]{90}{GDSW}}}
& $4\times 4\times 4$ & 279   & { 29.49} & { 33} & { 57.82} & { 41} & { 63.46} & { 42} & { 66.57} & { 44} & { 68.78} & { 45} \\
& & $6\times 6\times 6$ & 1115   & { 24.99} & { 32} & { 54.95} & { 42} & { 61.06} & { 43} & { 64.40} & { 45} & { 66.94} & { 46} \\
& & $8\times 8\times 8$ & 2863   & { 21.40} & { 31} & { 52.02} & { 42} & { 58.76} & { 44} & { 62.28} & { 45} & { 65.04} & { 47} \\
& & $12\times 12\times 12$ & 10439   & { 15.87} & { 28} & { 46.22} & { 42} & { 54.41} & { 44} & { 58.16} & { 46} & { 61.27} & { 47} \\
\cline{2-14}
& \multirow{4}{*}{\STAB{\rotatebox[origin=c]{90}{GDSW$^*$}}}
& $4\times 4\times 4$ & 171   & { 29.11} & { 31} & { 56.95} & { 40} & { 62.98} & { 43} & { 65.71} & { 44} & { 67.76} & { 46} \\
& & $6\times 6\times 6$ & 665   & { 26.37} & { 33} & { 56.59} & { 42} & { 63.35} & { 44} & { 66.13} & { 47} & { 68.42} & { 47} \\
& & $8\times 8\times 8$ & 1687   & { 23.42} & { 32} & { 54.85} & { 44} & { 62.37} & { 46} & { 65.37} & { 48} & { 67.87} & { 49} \\
& & $12\times 12\times 12$ & 6083   & { 17.71} & { 30} & { 49.39} & { 45} & { 58.68} & { 48} & { 62.13} & { 49} & { 65.12} & { 51} \\
\cline{2-14}
& \multirow{4}{*}{\STAB{\rotatebox[origin=c]{90}{RGDSW}}}
& $4\times 4\times 4$ & 27   & { 37.56} & { 31} & { 71.05} & { 40} & { 78.35} & { 42} & { 81.89} & { 45} & { 84.29} & { 46} \\
& & $6\times 6\times 6$ & 125   & { 33.02} & { 33} & { 69.20} & { 43} & { 76.60} & { 45} & { 80.26} & { 47} & { 82.98} & { 48} \\
& & $8\times 8\times 8$ & 343   & { 28.48} & { 34} & { 65.55} & { 45} & { 73.40} & { 48} & { 77.23} & { 50} & { 80.14} & { 51} \\
& & $12\times 12\times 12$ & 1331   & { 20.72} & { 32} & { 57.08} & { 46} & { 66.35} & { 49} & { 70.55} & { 51} & { 73.91} & { 52} \\
\hline
\end{tabular}
\end{adjustbox}
\end{table}
\begin{figure}[tbp]
\centering
\includegraphics[width=\textwidth]{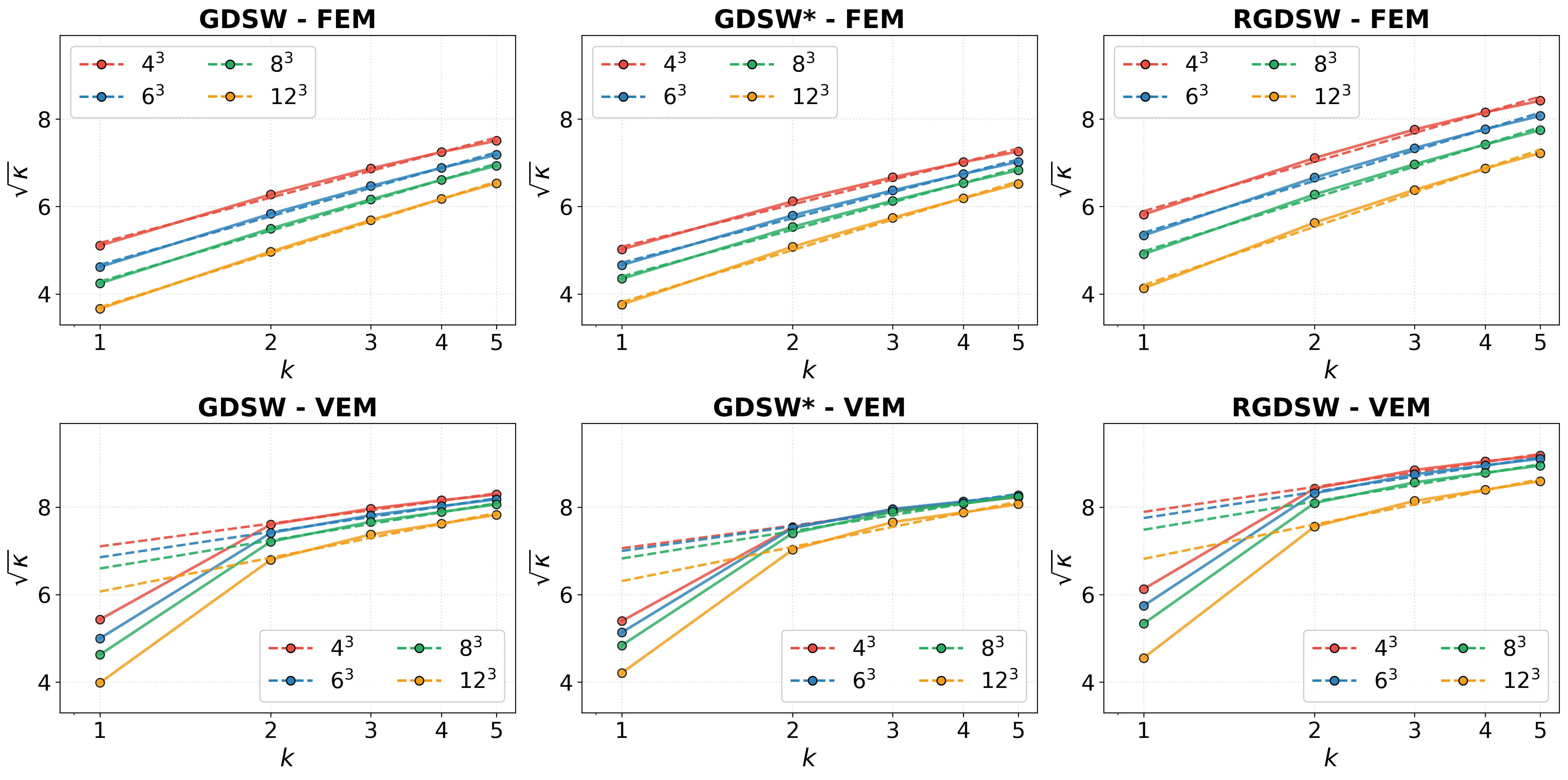}
\vspace{-7mm}
\caption{Square roots of the condition numbers from \mbox{\Cref{tab:3d_160_kappa}} for $\mathcal{Q}_k$-FEM and \ac{VEM} with GDSW, GDSW$^*$, and RGDSW. Dashed curves show least-squares fits against $1+\log(Hk^3/h)$.}
\label{fig:3d_table1_kappa}
\end{figure}
\Cref{tab:3d_160_kappa} shows condition numbers and iteration counts for \ac{FEM} and \ac{VEM} on a fixed $48\times48\times48$ mesh, varying the number of subdomains and the polynomial degree $k$.
For both \ac{FEM} and \ac{VEM}, increasing the number of subdomains does not lead to a substantial growth of the iteration counts.
For \ac{FEM}, the growth of the condition number as the polynomial degree increases from 1 to 5 is consistent with logarithmic-type behavior. \mbox{\Cref{fig:3d_table1_kappa}} illustrates this through plots of $\sqrt{\kappa}$ and least-squares fits against $1+\log(Hk^3/h)$.
In the absence of a theoretical condition number bound in three dimensions, this fit was chosen empirically.
For \ac{VEM}, the behavior is qualitatively similar, with the exception of a sharper growth from $k=1$ to $k=2$ followed by a smaller growth than for \ac{FEM}.
As expected, condition numbers are overall higher than those of their \ac{FEM} counterparts at equal degrees.
Among the three variants, GDSW and GDSW$^*$ yield the smallest condition numbers and are the most robust choices.
RGDSW, the cheapest option, shows a significant coarse dimension reduction, at the price of a modest increase in condition numbers and a few additional iterations.

\subsection{Parallel implementation strong scaling}
\label{sec:results_3d_strong}

\mbox{\Cref{tab:strong_fem_fair}} reports parallel strong scaling results for \ac{FEM} and \ac{VEM} with GDSW, GDSW$^*$, and RGDSW at degree $k=5$.
We fix the decomposition to 768 subdomains and double the number of computational nodes from 1 to 8, keeping 96 cores per node (96, 192, 384, and 768 total cores).
The starting configuration (1 node, 96 cores) uses the largest problem that fits on a single node with one terabyte of memory: $7\,312\,581$ \ac{DOFs} for \ac{FEM} and $3\,714\,501$ for \ac{VEM}.
We compute the parallel efficiency $E_{\rm p}$ on the solver time $T_{\rm setup}+T_{\rm sol}$, i.e., the part of the computation handled by the preconditioner.
For \ac{FEM}, all three variants exhibit good scalability, with $E_{\rm p}$ remaining at about $90\%$ at 8 nodes.
For \ac{VEM}, $E_{\rm p}$ stays above $91\%$ up to 4 nodes and ranges from $84\%$ to $87\%$ at 8 nodes, depending on the variant.
The assembly time, which is local to the subdomains, scales almost ideally for both discretizations.
The significantly higher \ac{VEM} assembly time is due to the element-level computations in \texttt{VEM++}, which have not been optimized; further reductions in \ac{VEM} assembly time would require optimizations in the \texttt{VEM++} library, which is outside the scope of this work.
Regarding the coarse setup, RGDSW consistently reduces $T_{\rm coarse}$ compared to GDSW across all node counts, while GDSW$^*$ is intermediate between the two, reflecting the relative sizes of their coarse spaces.

\begin{table}[tp]
\centering
\caption{Parallel strong scaling for \ac{FEM} and \ac{VEM} with GDSW, GDSW$^*$, and RGDSW for $k=5$ and 768 subdomains. We keep 96 cores per node and we increase the number of nodes. \ac{FEM}: \ac{DOFs} $7\,312\,581$, mesh $40\times40\times36$. \ac{VEM}: \ac{DOFs} $3\,714\,501$, mesh $40\times40\times36$. All times in seconds; $T_{\rm ass}$ denotes assembly time, $T_{\rm setup}$ preconditioner setup time, $T_{\rm sub}$ setup time of the local subdomain problems (factorizations), $T_{\rm coarse}$ setup time of the coarse problem (coarse basis, $A_0$, and its factorization), $T_{\rm sol}$ PCG solve time, and $E_{\rm p}$ parallel efficiency computed from $T_{\rm setup}+T_{\rm sol}$.}\label{tab:strong_fem_fair}
\begin{adjustbox}{width=\textwidth}
\begin{tabular}{|c|c|r|r||r|r|r|r|r|r|}
\cline{3-10}
\multicolumn{2}{c|}{} & Nodes & Cores & $T_{\rm ass}$ & $T_{\rm setup}$ & $T_{\rm sub}$ & $T_{\rm coarse}$ & $T_{\rm sol}$ & $E_{\rm p}$ \\
\cline{3-10}
\noalign{\vskip 0.6ex}
\cline{1-10}
\multirow{12}{*}{\STAB{\rotatebox[origin=c]{90}{FEM}}}
& \multirow{4}{*}{\STAB{\rotatebox[origin=c]{90}{GDSW}}}
& 1 & 96 & 5.20 & 34.36 & 23.15 & 11.21 & 18.24 & 1.000 \\
& & 2 & 192 & 2.68 & 17.90 & 12.15 & 5.76 & 9.31 & 0.967 \\
& & 4 & 384 & 1.45 & 9.18 & 6.18 & 3.00 & 4.80 & 0.941 \\
& & 8 & 768 & 0.88 & 4.77 & 3.17 & 1.60 & 2.57 & 0.896 \\
\cline{2-10}
& \multirow{4}{*}{\STAB{\rotatebox[origin=c]{90}{GDSW$^*$}}}
& 1 & 96 & 4.97 & 32.63 & 23.00 & 9.62 & 19.01 & 1.000 \\
& & 2 & 192 & 2.70 & 17.06 & 12.20 & 4.87 & 9.70 & 0.965 \\
& & 4 & 384 & 1.65 & 8.74 & 6.20 & 2.54 & 4.98 & 0.941 \\
& & 8 & 768 & 1.03 & 4.54 & 3.21 & 1.34 & 2.66 & 0.896 \\
\cline{2-10}
& \multirow{4}{*}{\STAB{\rotatebox[origin=c]{90}{RGDSW}}}
& 1 & 96 & 4.98 & 31.61 & 23.03 & 8.58 & 20.21 & 1.000 \\
& & 2 & 192 & 2.68 & 16.46 & 12.13 & 4.34 & 10.29 & 0.968 \\
& & 4 & 384 & 1.68 & 8.39 & 6.15 & 2.24 & 5.28 & 0.947 \\
& & 8 & 768 & 1.02 & 4.40 & 3.22 & 1.17 & 2.80 & 0.901 \\
\hline
\hline
\multirow{12}{*}{\STAB{\rotatebox[origin=c]{90}{VEM}}}
& \multirow{4}{*}{\STAB{\rotatebox[origin=c]{90}{GDSW}}}
& 1 & 96 & 256.12 & 14.22 & 9.49 & 4.73 & 11.27 & 1.000 \\
& & 2 & 192 & 127.74 & 7.50 & 4.99 & 2.50 & 5.81 & 0.958 \\
& & 4 & 384 & 65.33 & 3.88 & 2.48 & 1.40 & 3.03 & 0.922 \\
& & 8 & 768 & 32.74 & 2.10 & 1.33 & 0.77 & 1.67 & 0.845 \\
\cline{2-10}
& \multirow{4}{*}{\STAB{\rotatebox[origin=c]{90}{GDSW$^*$}}}
& 1 & 96 & 257.03 & 13.29 & 9.39 & 3.90 & 12.17 & 1.000 \\
& & 2 & 192 & 127.88 & 6.95 & 4.92 & 2.03 & 6.24 & 0.965 \\
& & 4 & 384 & 63.99 & 3.57 & 2.49 & 1.08 & 3.25 & 0.934 \\
& & 8 & 768 & 34.83 & 1.92 & 1.32 & 0.60 & 1.78 & 0.862 \\
\cline{2-10}
& \multirow{4}{*}{\STAB{\rotatebox[origin=c]{90}{RGDSW}}}
& 1 & 96 & 255.69 & 12.84 & 9.31 & 3.52 & 12.90 & 1.000 \\
& & 2 & 192 & 128.06 & 6.71 & 4.89 & 1.82 & 6.60 & 0.967 \\
& & 4 & 384 & 66.59 & 3.44 & 2.48 & 0.97 & 3.44 & 0.935 \\
& & 8 & 768 & 32.60 & 1.85 & 1.32 & 0.53 & 1.87 & 0.866 \\
\hline
\end{tabular}
\end{adjustbox}
\end{table}

\subsection{A weak scaling test}
\label{sec:results_3d_scal}

In this test, we fix the ratio $H/h=4$, i.e., $4^3$ elements per subdomain, and consider $4^3$, $8^3$, $12^3$, and $16^3$ subdomains on the unit cube, for degree $k=5$.
\Cref{tab:scalability} shows that the condition numbers and the iteration counts still grow with the number of subdomains, but only moderately. The setup time of the local problems stays almost constant, and the time per PCG iteration grows only slightly, so that the increase of $T_{\rm sol}$ is mainly due to the iteration counts. As expected, the setup time of the coarse problem grows with the coarse space dimension. For GDSW with $4\,096$ subdomains, it reaches the local setup time for \ac{VEM} and about two thirds of it for \ac{FEM}, whereas for RGDSW it remains almost constant, at the price of a few additional iterations. Since the coarse problem is solved directly, a constant time to solution cannot be expected.

\begin{table}[tbp]

\centering
\caption{Numerical scalability for \ac{FEM} and \ac{VEM} with GDSW, GDSW$^*$, and RGDSW on the unit cube for $k=5$ and $H/h=4$, with one subdomain per core. $T_{\rm sub}$ denotes the setup time of the local subdomain problems, $T_{\rm coarse}$ the setup time of the coarse problem, and $T_{\rm sol}$ the PCG solve time; all times in seconds.}\label{tab:scalability}
\begin{adjustbox}{width=\textwidth}
\begin{tabular}{|c|c|r|r|r|r||r|r||r|r|r|}
\cline{3-11}
\multicolumn{2}{c|}{} & nSub & Cores & \ac{DOFs} & Coarse & $\kappa$ & it & $T_{\rm sub}$ & $T_{\rm coarse}$ & $T_{\rm sol}$ \\
\cline{3-11}
\noalign{\vskip 0.6ex}
\cline{1-11}
\multirow{12}{*}{\STAB{\rotatebox[origin=c]{90}{FEM}}}
& \multirow{4}{*}{\STAB{\rotatebox[origin=c]{90}{GDSW}}}
& $4\times 4\times 4$ & 64 & 531\,441 & 279 & 24.12 & 28 & 2.49 & 1.28 & 1.41 \\
& & $8\times 8\times 8$ & 512 & 4\,173\,281 & 2863 & 37.75 & 36 & 2.89 & 1.40 & 2.11 \\
& & $12\times 12\times 12$ & 1728 & 13\,997\,521 & 10439 & 42.70 & 40 & 3.05 & 1.62 & 2.48 \\
& & $16\times 16\times 16$ & 4096 & 33\,076\,161 & 25695 & 44.75 & 45 & 3.06 & 2.00 & 2.95 \\
\cline{2-11}
& \multirow{4}{*}{\STAB{\rotatebox[origin=c]{90}{GDSW$^*$}}}
& $4\times 4\times 4$ & 64 & 531\,441 & 171 & 23.50 & 28 & 2.52 & 1.16 & 1.41 \\
& & $8\times 8\times 8$ & 512 & 4\,173\,281 & 1687 & 37.20 & 36 & 2.83 & 1.23 & 2.08 \\
& & $12\times 12\times 12$ & 1728 & 13\,997\,521 & 6083 & 42.51 & 42 & 3.02 & 1.30 & 2.57 \\
& & $16\times 16\times 16$ & 4096 & 33\,076\,161 & 14895 & 44.74 & 45 & 3.07 & 1.46 & 2.84 \\
\cline{2-11}
& \multirow{4}{*}{\STAB{\rotatebox[origin=c]{90}{RGDSW}}}
& $4\times 4\times 4$ & 64 & 531\,441 & 27 & 27.21 & 26 & 2.52 & 1.05 & 1.30 \\
& & $8\times 8\times 8$ & 512 & 4\,173\,281 & 343 & 45.91 & 38 & 2.89 & 1.10 & 2.18 \\
& & $12\times 12\times 12$ & 1728 & 13\,997\,521 & 1331 & 52.17 & 45 & 3.08 & 1.12 & 2.70 \\
& & $16\times 16\times 16$ & 4096 & 33\,076\,161 & 3375 & 54.87 & 51 & 3.12 & 1.19 & 3.12 \\
\hline
\hline
\multirow{12}{*}{\STAB{\rotatebox[origin=c]{90}{VEM}}}
& \multirow{4}{*}{\STAB{\rotatebox[origin=c]{90}{GDSW}}}
& $4\times 4\times 4$ & 64 & 272\,881 & 279 & 29.20 & 30 & 1.06 & 0.65 & 0.78 \\
& & $8\times 8\times 8$ & 512 & 2\,123\,233 & 2863 & 51.64 & 39 & 1.15 & 0.60 & 1.10 \\
& & $12\times 12\times 12$ & 1728 & 7\,099\,345 & 10439 & 61.27 & 47 & 1.20 & 0.80 & 1.46 \\
& & $16\times 16\times 16$ & 4096 & 16\,749\,505 & 25695 & 65.57 & 52 & 1.23 & 1.20 & 1.74 \\
\cline{2-11}
& \multirow{4}{*}{\STAB{\rotatebox[origin=c]{90}{GDSW$^*$}}}
& $4\times 4\times 4$ & 64 & 272\,881 & 171 & 29.26 & 31 & 1.04 & 0.53 & 0.79 \\
& & $8\times 8\times 8$ & 512 & 2\,123\,233 & 1687 & 54.41 & 42 & 1.15 & 0.48 & 1.17 \\
& & $12\times 12\times 12$ & 1728 & 7\,099\,345 & 6083 & 65.12 & 51 & 1.20 & 0.54 & 1.52 \\
& & $16\times 16\times 16$ & 4096 & 16\,749\,505 & 14895 & 69.72 & 57 & 1.23 & 0.69 & 1.82 \\
\cline{2-11}
& \multirow{4}{*}{\STAB{\rotatebox[origin=c]{90}{RGDSW}}}
& $4\times 4\times 4$ & 64 & 272\,881 & 27 & 31.97 & 29 & 1.06 & 0.46 & 0.75 \\
& & $8\times 8\times 8$ & 512 & 2\,123\,233 & 343 & 61.86 & 43 & 1.17 & 0.42 & 1.18 \\
& & $12\times 12\times 12$ & 1728 & 7\,099\,345 & 1331 & 73.91 & 52 & 1.21 & 0.44 & 1.52 \\
& & $16\times 16\times 16$ & 4096 & 16\,749\,505 & 3375 & 78.94 & 60 & 1.30 & 0.48 & 1.81 \\
\hline
\end{tabular}
\end{adjustbox}
\end{table}

\subsection{Necessity of the moment-based face evaluation strategy for high-order \ac{VEM}}
\label{sec:results_3d_strip}

\begin{figure}[tbp]
\centering
\includegraphics[trim=0 200 0 250, clip,width=\textwidth]{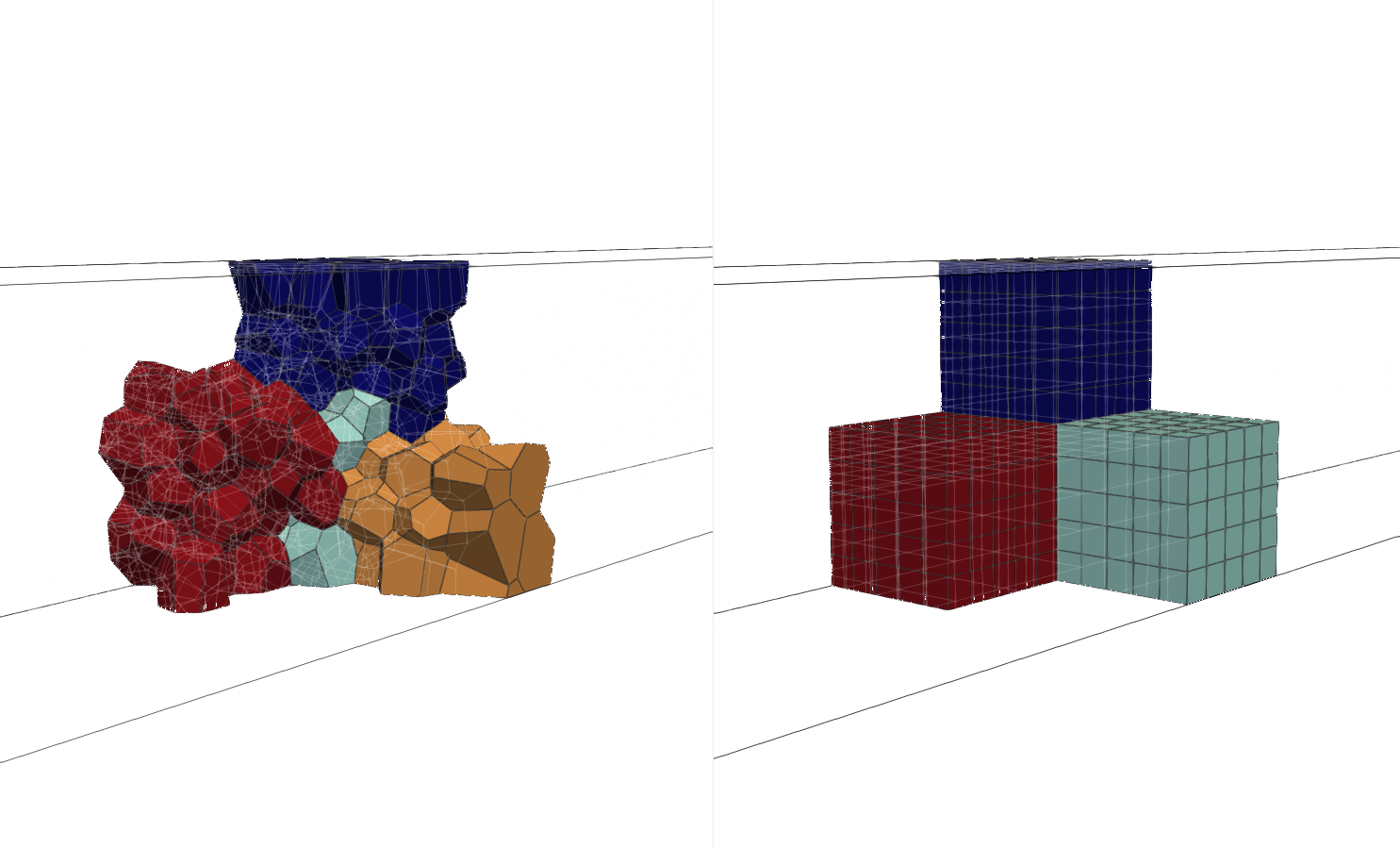}
\caption{Local subdomains of $\Omega=(0,64)\times(0,1)\times(0,1)$. Left: CVT polyhedral mesh used for the \ac{VEM} discretization, and the corresponding domain decomposition, in \cref{sec:results_3d_strip}. Right: CUBE mesh and corresponding subdomain decomposition used in \cref{sec:results_3d_strip,sec:results_3d_tol,sec:results_channel}.}
\label{fig:slab_domain}
\end{figure}

In this experiment, we demonstrate the necessity of the moment-based face evaluation introduced in \Cref{sec:gdsw_high_order} for high-order \ac{VEM} by comparing two families of coarse spaces.
The standard variants, denoted GDSW$^{\rm std}$, GDSW$^{*,\rm std}$, and RGDSW$^{\rm std}$, apply the \ac{FEM}-style face evaluation, i.e., they assign the value $1$ to all interface \ac{DOFs}. This is correct for nodal \ac{FEM} and for \ac{VEM} up to degree $k=2$, but fails for high-order \ac{VEM} with $k\ge3$, where the vector of all ones no longer represents the constant function.
The corrected variants, denoted GDSW, GDSW$^*$, and RGDSW, use the moment-based face evaluation of \Cref{sec:gdsw_high_order}.
We consider two mesh types on the slab domain (cf.\ \cref{fig:slab_domain}): CUBE ($128\times2\times2$ hexahedral decomposition with $5^3$ elements per subdomain) and CVT ($16\,000$ polyhedral elements generated with 20 Lloyd iterations~\cite{du_centroidal_1999} and partitioned by METIS~\cite{karypis_metis_1998}), both with 512 subdomains.

\begin{table}[tbp]
\centering
\caption{Condition numbers and PCG iteration counts for \ac{VEM} on the slab domain with 512 subdomains and CUBE/CVT meshes, for polynomial degrees $k=1,\ldots,5$. The standard coarse spaces (superscript ${\rm std}$) use \ac{FEM}-style face evaluation, while the corrected variants use the moment-based face evaluation of \cref{sec:gdsw_high_order}.}
\label{tab:degree}
\begin{adjustbox}{width=\textwidth}
\begin{tabular}{|c|l||r||r|r||r|r||r|r||r|r||r|r|}
\cline{3-13}
\multicolumn{2}{c|}{} & Coarse & \multicolumn{2}{c||}{$k=1$} & \multicolumn{2}{c||}{$k=2$} & \multicolumn{2}{c||}{$k=3$} & \multicolumn{2}{c||}{$k=4$} & \multicolumn{2}{c|}{$k=5$} \\
\cline{3-13}
\multicolumn{2}{c|}{} & & $\kappa$ & it & $\kappa$ & it & $\kappa$ & it & $\kappa$ & it & $\kappa$ & it \\
\cline{3-13}
\noalign{\vskip 0.6ex}
\cline{1-13}
\hline
\multirow{7}{*}{\STAB{\rotatebox[origin=c]{90}{CUBE}}}
& ASM (1-lv) & \multicolumn{1}{c||}{--} & $1.2$e4 & $147$ & $2.2$e4 & $205$ & $2.5$e4 & $221$ & $2.6$e4 & $241$ & $2.6$e4 & $250$ \\
& GDSW$^{\rm std}$ & $1783$ & $16.2$ & $30$ & $36.3$ & $43$ & $2.8$e4 & $269$ & $2.9$e4 & $283$ & $3.0$e4 & $290$ \\
& GDSW$^{*,\rm std}$ & $1147$ & $15.7$ & $29$ & $35.1$ & $43$ & $2.6$e4 & $259$ & $2.7$e4 & $275$ & $2.7$e4 & $281$ \\
& RGDSW$^{\rm std}$ & $127$ & $19.7$ & $28$ & $43.2$ & $42$ & $2.5$e4 & $249$ & $2.6$e4 & $264$ & $2.7$e4 & $272$ \\
& GDSW & $1783$ & $16.2$ & $30$ & $36.3$ & $43$ & $40.3$ & $46$ & $42.5$ & $49$ & $44.4$ & $50$ \\
& GDSW$^*$ & $1147$ & $15.7$ & $29$ & $35.1$ & $43$ & $39.1$ & $45$ & $41.2$ & $48$ & $43.0$ & $49$ \\
& RGDSW & $127$ & $19.7$ & $28$ & $43.2$ & $42$ & $48.1$ & $45$ & $50.7$ & $47$ & $52.8$ & $49$ \\
\hline
\multirow{7}{*}{\STAB{\rotatebox[origin=c]{90}{CVT}}}
& ASM (1-lv) & \multicolumn{1}{c||}{--} & $3.8$e4 & $326$ & $7.0$e4 & $450$ & $8.1$e4 & $485$ & $9.3$e4 & $518$ & $9.9$e4 & $529$ \\
& GDSW$^{\rm std}$ & $4902$ & $26.5$ & $44$ & $49.8$ & $60$ & $8.6$e4 & $496$ & $9.3$e4 & $516$ & $9.8$e4 & $529$ \\
& GDSW$^{*,\rm std}$ & $3067$ & $28.2$ & $45$ & $54.7$ & $63$ & $8.6$e4 & $493$ & $9.3$e4 & $515$ & $9.8$e4 & $529$ \\
& RGDSW$^{\rm std}$ & $1468$ & $31.3$ & $47$ & $59.6$ & $65$ & $8.6$e4 & $493$ & $9.3$e4 & $515$ & $9.8$e4 & $529$ \\
& GDSW & $4902$ & $26.5$ & $44$ & $49.8$ & $60$ & $52.8$ & $62$ & $55.8$ & $63$ & $59.7$ & $66$ \\
& GDSW$^*$ & $3067$ & $28.2$ & $45$ & $54.7$ & $63$ & $57.1$ & $64$ & $61.3$ & $66$ & $65.3$ & $69$ \\
& RGDSW & $1468$ & $31.3$ & $47$ & $59.6$ & $65$ & $63.0$ & $67$ & $68.8$ & $70$ & $72.8$ & $72$ \\
\hline
\end{tabular}
\end{adjustbox}
\end{table}
The results in~\cref{tab:degree} confirm the necessity of the moment-based construction of \Cref{sec:gdsw_high_order}. At $k=1,2$ all variants coincide, as expected. For $k\ge 3$, the standard variants GDSW$^{\rm std}$, GDSW$^{*,\rm std}$, and RGDSW$^{\rm std}$ fail on both mesh types, with condition numbers matching those of the one-level solver. 
Without the moment-based face evaluation, the coarse space does not contain the constant functions, i.e., the kernel of the operator, and the coarse correction does not provide the desired effect. Introducing the moment-based face evaluation restores the quality of the coarse space, and all three variants maintain bounded condition numbers across all degrees and mesh types.

\subsection{Effect of coarse basis truncation}
\label{sec:results_3d_tol}

In this test, we use the cube mesh and a regular decomposition into cube subdomains and elements (cf.\ \cref{fig:slab_domain}, right), with $512$ subdomains and $512$ MPI ranks.
In \mbox{\Cref{tab:tol}}, we study the effect of sparsifying the coarse basis matrix $\Phi$ by zeroing entries below the threshold $\phi_{\rm tol}$, which reduces the cost of the coarse setup at the price of a less accurate coarse space.
For \ac{FEM}, none of the variants is affected up to $\phi_{\rm tol}=10^{-4}$, while $\phi_{\rm tol}=10^{-2}$ degrades GDSW and GDSW$^*$ for all degrees.
For \ac{VEM}, the moment-based face values introduce small entries in $\Phi$, which are removed already by moderate thresholds. These entries are small but not negligible: the rows of $\Phi$ associated with moment \ac{DOFs} are scaled differently from the nodal rows, so that the same absolute threshold is not comparable between \ac{FEM} and \ac{VEM}. GDSW and GDSW$^*$ are unaffected by $\phi_{\rm tol}=10^{-6}$ up to $k=4$, while $\phi_{\rm tol}=10^{-4}$ degrades them for $k\ge 3$.
RGDSW is the least sensitive variant since its basis functions group vertex, edge, and face contributions; therefore, the entries of $\Phi$ are larger and fewer of them are filtered out.

\begin{table}[tbp]
\centering
\caption{{Condition numbers and PCG iteration counts for \ac{FEM} and \ac{VEM} with all GDSW variants at varying truncation threshold}
  $\phi_{\rm tol}$. The test configuration is the same as in \mbox{\Cref{tab:degree}}, with a CUBE mesh and a regular decomposition into $512$ cubic subdomains.}
\label{tab:tol}
\begin{adjustbox}{width=\textwidth}
\begin{tabular}{|c|c|c||r|r||r|r||r|r||r|r||r|r|}
\cline{4-13}
\multicolumn{3}{c|}{} & \multicolumn{2}{c||}{$k=1$} & \multicolumn{2}{c||}{$k=2$} & \multicolumn{2}{c||}{$k=3$} & \multicolumn{2}{c||}{$k=4$} & \multicolumn{2}{c|}{$k=5$} \\
\cline{3-13}
\multicolumn{2}{c|}{} & ${\phi_{\rm tol}}$ & $\kappa$ & it & $\kappa$ & it & $\kappa$ & it & $\kappa$ & it & $\kappa$ & it \\
\cline{3-13}
\noalign{\vskip 0.6ex}
\cline{1-13}
\multirow{12}{*}{\STAB{\rotatebox[origin=c]{90}{FEM}}}
& \multirow{4}{*}{GDSW} & $0$ & $13.8$ & $27$ & $21.1$ & $33$ & $25.8$ & $37$ & $29.3$ & $39$ & $32.1$ & $40$ \\
& & $10^{-6}$ & $13.8$ & $27$ & $21.1$ & $33$ & $25.8$ & $37$ & $29.3$ & $39$ & $32.1$ & $40$ \\
& & $10^{-4}$ & $13.8$ & $27$ & $21.1$ & $33$ & $25.8$ & $37$ & $29.3$ & $39$ & $32.1$ & $40$ \\
& & $10^{-2}$ & $39.3$ & $35$ & $28.2$ & $36$ & $52.4$ & $43$ & $96.9$ & $50$ & $169.8$ & $57$ \\
\cline{2-13}
& \multirow{4}{*}{GDSW$^*$} & $0$ & $13.2$ & $27$ & $19.9$ & $33$ & $24.2$ & $37$ & $27.5$ & $38$ & $30.1$ & $40$ \\
& & $10^{-6}$ & $13.2$ & $27$ & $19.9$ & $33$ & $24.2$ & $37$ & $27.5$ & $38$ & $30.1$ & $40$ \\
& & $10^{-4}$ & $13.2$ & $27$ & $19.9$ & $33$ & $24.2$ & $37$ & $27.5$ & $38$ & $30.1$ & $40$ \\
& & $10^{-2}$ & $18.4$ & $29$ & $25.9$ & $35$ & $48.3$ & $43$ & $89.3$ & $49$ & $156.5$ & $57$ \\
\cline{2-13}
& \multirow{4}{*}{RGDSW} & $0$ & $17.1$ & $27$ & $26.3$ & $33$ & $32.1$ & $35$ & $36.3$ & $39$ & $39.6$ & $40$ \\
& & $10^{-6}$ & $17.1$ & $27$ & $26.3$ & $33$ & $32.1$ & $35$ & $36.3$ & $39$ & $39.6$ & $40$ \\
& & $10^{-4}$ & $17.1$ & $27$ & $26.3$ & $33$ & $32.1$ & $35$ & $36.3$ & $39$ & $39.6$ & $40$ \\
& & $10^{-2}$ & $17.1$ & $27$ & $26.3$ & $33$ & $32.1$ & $35$ & $36.3$ & $39$ & $39.6$ & $40$ \\
\hline\hline
\multirow{12}{*}{\STAB{\rotatebox[origin=c]{90}{VEM}}}
& \multirow{4}{*}{GDSW} & $0$ & $16.2$ & $30$ & $36.3$ & $43$ & $40.3$ & $46$ & $42.5$ & $49$ & $44.4$ & $50$ \\
& & $10^{-6}$ & $16.2$ & $30$ & $36.3$ & $43$ & $40.3$ & $46$ & $42.7$ & $49$ & $72.2$ & $56$ \\
& & $10^{-4}$ & $16.2$ & $30$ & $36.3$ & $43$ & $65.4$ & $52$ & $1986$ & $116$ & $2.7$e4 & $280$ \\
& & $10^{-2}$ & $131.3$ & $43$ & $806.4$ & $86$ & $1.9$e4 & $237$ & $2.9$e4 & $285$ & $3.0$e4 & $287$ \\
\cline{2-13}
& \multirow{4}{*}{GDSW$^*$} & $0$ & $15.7$ & $29$ & $35.1$ & $43$ & $39.1$ & $45$ & $41.2$ & $48$ & $43.0$ & $49$ \\
& & $10^{-6}$ & $15.7$ & $29$ & $35.1$ & $43$ & $39.1$ & $45$ & $41.2$ & $48$ & $60.9$ & $55$ \\
& & $10^{-4}$ & $15.7$ & $29$ & $35.1$ & $43$ & $47.4$ & $48$ & $1415$ & $107$ & $2.5$e4 & $271$ \\
& & $10^{-2}$ & $21.5$ & $31$ & $394.6$ & $73$ & $1.8$e4 & $230$ & $2.7$e4 & $279$ & $2.8$e4 & $283$ \\
\cline{2-13}
& \multirow{4}{*}{RGDSW} & $0$ & $19.7$ & $28$ & $43.2$ & $42$ & $48.1$ & $45$ & $50.7$ & $47$ & $52.8$ & $49$ \\
& & $10^{-6}$ & $19.7$ & $28$ & $43.2$ & $42$ & $48.1$ & $45$ & $50.7$ & $47$ & $52.8$ & $49$ \\
& & $10^{-4}$ & $19.7$ & $28$ & $43.2$ & $42$ & $48.1$ & $45$ & $51.1$ & $47$ & $116.5$ & $73$ \\
& & $10^{-2}$ & $19.7$ & $28$ & $43.2$ & $42$ & $134.0$ & $77$ & $2.6$e4 & $257$ & $2.7$e4 & $269$ \\
\hline
\end{tabular}
\end{adjustbox}
\end{table}

\subsection{Heterogeneous coefficients: Four-channel problem}
\label{sec:results_channel}

We now consider the scalar diffusion problem \eqref{eq:diffusion} with a heterogeneous coefficient function $\rho$.

\begin{figure}[tbp]
\centering
\includegraphics[trim=100 325 650 325, clip, width=0.9\textwidth]{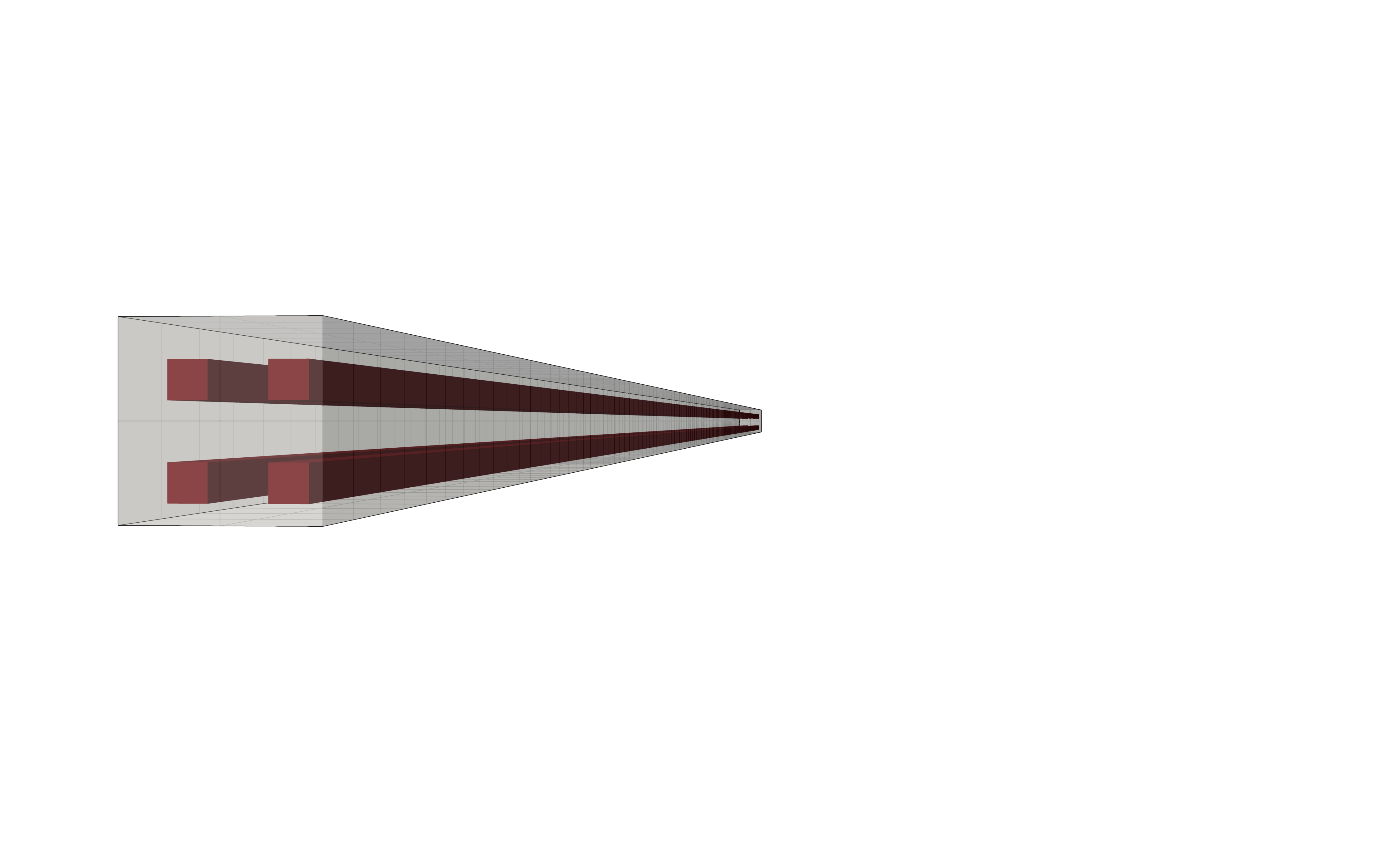}
\caption{Coefficient field and subdomain decomposition for the four-channel problem. The grid lines show the subdomain boundaries, the four channels (red) have $\rho_{\text{max}}=1000$, while the background (grey) has $\rho_{\text{min}}=1$.}
\label{fig:channels}
\end{figure}

In the experiment below, the permeability takes different values inside and outside the channels, with the channel boundaries crossing the subdomain interfaces; see \mbox{\Cref{fig:channels}}.
The domain is the slab $64\times 1\times 1$ with four channels oriented along the $x$-axis, each with a square cross-section of side $0.2$, occupying $y,z\in[0.1,0.3]\cup[0.6,0.8]$, i.e., two layers of elements in each direction, running from $x=0.1$ to $x=63.9$ (cf.\ \cref{fig:channels}).
The permeability function is $\rho=\rho_{\text{max}}\,I$ inside the channels and $\rho=\rho_{\text{min}}\,I$ outside, with $\rho_{\text{max}}=1000$ and $\rho_{\text{min}}=1$.
Here, we again use the CUBE mesh and the regular decomposition into cubic subdomains and elements (cf.\ \cref{fig:slab_domain}, right).

\begin{table}[tbp]
\centering
\caption{Condition numbers and PCG iteration counts for \ac{FEM} and \ac{VEM} on the four-channel heterogeneous diffusion problem ($\rho_{\text{max}}/\rho_{\text{min}}=1000$) on the slab domain with a CUBE mesh and a regular cubic decomposition into 512 subdomains, for polynomial degrees $k=1,\ldots,5$.}
\label{tab:channel}
\begin{adjustbox}{width=\textwidth}
\begin{tabular}{|c|l||r||r|r||r|r||r|r||r|r||r|r|}
\cline{3-13}
\multicolumn{2}{c|}{} & Coarse & \multicolumn{2}{c||}{$k=1$} & \multicolumn{2}{c||}{$k=2$} & \multicolumn{2}{c||}{$k=3$} & \multicolumn{2}{c||}{$k=4$} & \multicolumn{2}{c|}{$k=5$} \\
\cline{3-13}
\multicolumn{2}{c|}{} & & $\kappa$ & it & $\kappa$ & it & $\kappa$ & it & $\kappa$ & it & $\kappa$ & it \\
\cline{3-13}
\noalign{\vskip 0.6ex}
\cline{1-13}
\multirow{4}{*}{\STAB{\rotatebox[origin=c]{90}{FEM}}}
& ASM (1-lv) & \multicolumn{1}{c||}{--} & $2.0$e4 & $240$ & $2.8$e4 & $280$ & $3.0$e4 & $326$ & $3.0$e4 & $293$ & $3.6$e4 & $304$ \\
& GDSW  & $1783$ & $9.8$ & $24$ & $16.3$ & $27$ & $27.3$ & $37$ & $26.7$ & $36$ & $21.4$ & $34$ \\
& GDSW$^*$ & $1147$ & $11.1$ & $25$ & $20.2$ & $29$ & $28.5$ & $38$ & $28.1$ & $38$ & $22.3$ & $34$ \\
& RGDSW & $127$ & $1.1$e3 & $161$ & $2.2$e3 & $160$ & $1.9$e3 & $217$ & $1.3$e3 & $167$ & $1.7$e3 & $174$ \\
\hline\hline
\multirow{4}{*}{\STAB{\rotatebox[origin=c]{90}{VEM}}}
& ASM(1-lv) & \multicolumn{1}{c||}{--} & $1.9$e4 & $234$ & $3.2$e4 & $325$ & $3.7$e4 & $352$ & $3.8$e4 & $368$ & $3.9$e4 & $380$ \\
& GDSW & $1783$ & $12.6$ & $25$ & $33.1$ & $38$ & $37.6$ & $41$ & $40.8$ & $42$ & $43.7$ & $43$ \\
& GDSW$^*$ & $1147$ & $16.5$ & $28$ & $41.8$ & $40$ & $46.5$ & $43$ & $50.3$ & $44$ & $53.7$ & $47$ \\
& RGDSW & $127$ & $1.4$e3 & $132$ & $2.5$e3 & $178$ & $2.8$e3 & $199$ & $2.9$e3 & $210$ & $3.0$e3 & $214$ \\
\hline
\end{tabular}
\end{adjustbox}
\end{table}

\Cref{tab:channel} reports condition numbers and PCG iteration counts for $k=1,\ldots,5$.
For both \ac{FEM} and \ac{VEM}, GDSW, and GDSW$^*$ maintain bounded condition numbers for all degrees. For RGDSW, no robustness can be expected in this configuration. This is because that each interface plane is crossed by four channels along the $x$-direction, one per subdomain face, while RGDSW provides a single coarse function per interface plane, associated with its central vertex. GDSW and GDSW$^*$ have one coarse function per face and thus per channel. Accordingly, RGDSW is not robust for any polynomial degree, neither for \ac{FEM} nor for \ac{VEM}. GDSW and GDSW$^*$ remain robust, with condition numbers between $10$ and $30$ for \ac{FEM} and between $13$ and $54$ for \ac{VEM}.
Therefore, for the considered four-channel coefficient function, GDSW and GDSW$^*$ are the appropriate choices; RGDSW should be restricted to configurations where the coefficient is constant per subdomain or where the coefficient jumps align with subdomain boundaries.
In general, for arbitrary coefficient distributions with high permeability contrast ($\max_\Omega \rho / \min_\Omega \rho \gg 1$), robustness can only be guaranteed by adaptive coarse space strategies such as \ac{AGDSW}~\cite{heinlein_adaptive_2018,heinlein_adaptive_reduced_2022}, which enrich the coarse space with subdomain-local eigenvectors.

\subsection{Necessity of the moment-based face evaluation strategy for linear elasticity \ac{VEM}}
\label{sec:results_elasticity}

\begin{figure}[tbp]
\centering
\includegraphics[trim=-100 290 150 280, clip, width=\textwidth]{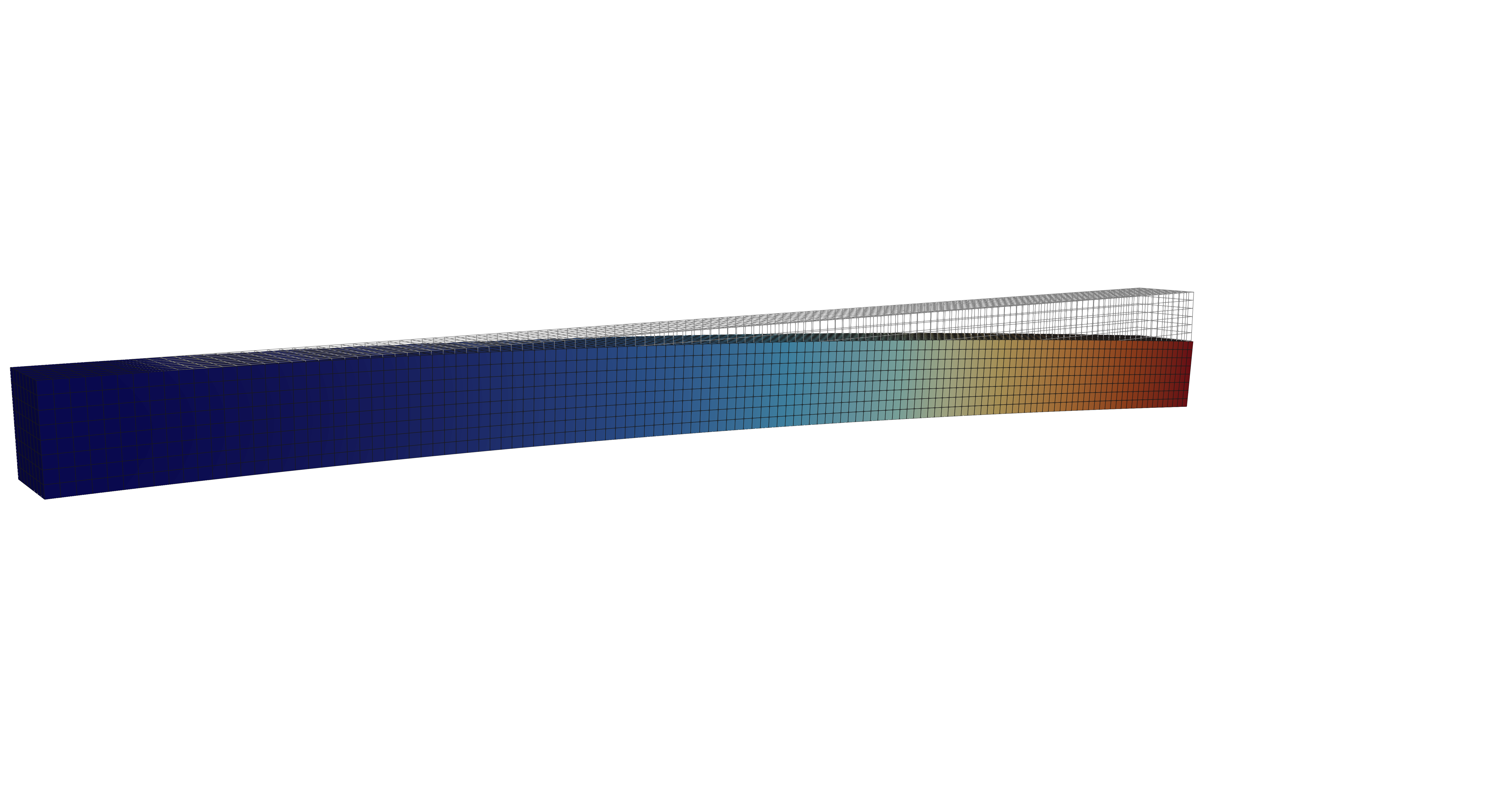}
\caption{Linear elasticity test on the slab domain: deformed configuration and displacement magnitude.}
\label{fig:lin_el_test}
\end{figure}

We extend the correctness test of \mbox{\Cref{sec:results_3d_strip}} to the vector-valued \ac{VEM} discretization of the linear elasticity problem introduced in \mbox{\Cref{sec:discretizations}}, using the componentwise moment-based face evaluation~\eqref{eq:gdsw_vem_face_el} described in \mbox{\Cref{sec:gdsw_high_order}}.
As in the scalar case, we compare the standard variants (GDSW$^{\rm std}$, GDSW$^{*,\rm std}$, and RGDSW$^{\rm std}$), which assign the value $1$ to all face \ac{DOFs}, against the corrected variants using the moment-based assignment.
The domain is the slab $\Omega=( 0,32)\times(0,1)\times(0,1)$ and we impose homogeneous Dirichlet boundary conditions $\mathbf{u}=\mathbf{0}$ on the face $\{x=0\}$, while Neumann boundary conditions are imposed on the remaining faces.
We apply a uniform body force $\mathbf{f}=(0,0,-1)^T$ to simulate gravity (cf.\ \Cref{fig:lin_el_test}).
The material parameters are Young's modulus $E=210\,000$ and Poisson's ratio $\nu=0.3$, corresponding to the Lam\'e parameters $\lambda = E\nu/((1+\nu)(1-2\nu))$ and $\mu = E/(2(1+\nu))$.
We consider two mesh configurations: a CUBE mesh with $256$ subdomains ($64\times2\times2$ decomposition) and $4^3$ elements per subdomain, and a CVT polyhedral mesh with $8\,000$ elements on the same domain with $256$ subdomains partitioned by METIS.

\begin{table}[tbp]
\centering
\caption{Condition numbers and PCG iteration counts for the linear elasticity problem on a CUBE mesh with $256$ subdomains ($64\times2\times2$), $4^3$ elements per subdomain, for $k=1,\ldots,5$.}
\label{tab:elasticity_degree}
\begin{adjustbox}{width=\textwidth}
\begin{tabular}{|c|l||r||r|r||r|r||r|r||r|r||r|r|}
\cline{3-13}
\multicolumn{2}{c|}{} & Coarse & \multicolumn{2}{c||}{$k=1$} & \multicolumn{2}{c||}{$k=2$} & \multicolumn{2}{c||}{$k=3$} & \multicolumn{2}{c||}{$k=4$} & \multicolumn{2}{c|}{$k=5$} \\
\cline{3-13}
\multicolumn{2}{c|}{} & & $\kappa$ & it & $\kappa$ & it & $\kappa$ & it & $\kappa$ & it & $\kappa$ & it \\
\cline{3-13}
\noalign{\vskip 0.6ex}
\cline{1-13}
\multirow{4}{*}{\STAB{\rotatebox[origin=c]{90}{FEM}}}
& ASM (1-lv) & \multicolumn{1}{c||}{--} & $4.8$e5 & $597$ & $5.7$e5 & $820$ & \multicolumn{2}{c||}{N.C.} & \multicolumn{2}{c||}{N.C.} & \multicolumn{2}{c|}{N.C.} \\
& GDSW  & $4817$ & $11.6$ & $27$ & $14.0$ & $30$ & $15.6$ & $31$ & $16.8$ & $33$ & $17.8$ & $34$ \\
& GDSW$^*$ & $3426$ & $11.4$ & $26$ & $14.1$ & $30$ & $15.7$ & $32$ & $16.9$ & $34$ & $17.8$ & $35$ \\
& RGDSW & $378$ & $15.0$ & $31$ & $21.9$ & $36$ & $26.5$ & $40$ & $29.7$ & $43$ & $32.1$ & $44$ \\
\hline
\hline
\multirow{7}{*}{\STAB{\rotatebox[origin=c]{90}{VEM}}}
& ASM (1-lv) & \multicolumn{1}{c||}{--} & $4.8$e5 & $646$ & $5.9$e5 & $992$ & \multicolumn{2}{c||}{N.C.} & \multicolumn{2}{c||}{N.C.} & \multicolumn{2}{c|}{N.C.} \\
& GDSW$^{\rm std}$ & $4817$ & $11.2$ & $26$ & $13.8$ & $30$ & \multicolumn{2}{c||}{N.C.} & \multicolumn{2}{c||}{N.C.} & \multicolumn{2}{c|}{N.C.} \\
& GDSW$^{*,\rm std}$ & $3426$ & $10.9$ & $26$ & $13.1$ & $30$ & \multicolumn{2}{c||}{N.C.} & \multicolumn{2}{c||}{N.C.} & \multicolumn{2}{c|}{N.C.} \\
& RGDSW$^{\rm std}$ & $378$ & $14.6$ & $31$ & $24.7$ & $38$ & \multicolumn{2}{c||}{N.C.} & \multicolumn{2}{c||}{N.C.} & \multicolumn{2}{c|}{N.C.} \\
& GDSW & $4817$ & $11.2$ & $26$ & $13.8$ & $30$ & $16.2$ & $33$ & $17.7$ & $35$ & $19.0$ & $35$ \\
& GDSW$^*$ & $3426$ & $10.9$ & $26$ & $13.1$ & $30$ & $15.3$ & $32$ & $16.8$ & $34$ & $18.1$ & $35$ \\
& RGDSW & $378$ & $14.6$ & $31$ & $24.7$ & $38$ & $30.7$ & $44$ & $33.9$ & $47$ & $36.6$ & $49$ \\
\hline
\end{tabular}
\end{adjustbox}
\end{table}

The results in \cref{tab:elasticity_degree} confirm that the moment-based face evaluation is likewise necessary for the \ac{VEM} discretization of linear elasticity.
As in the scalar diffusion case, at $k=1,2$ all variants coincide.
For $k\ge 3$, the standard variants GDSW$^{\rm std}$, GDSW$^{*,\rm std}$, and RGDSW$^{\rm std}$ do not converge, matching the one-level ASM performance.
The corrected variants using the moment-based assignment~\eqref{eq:gdsw_vem_face_el} maintain bounded condition numbers across all tested degrees, confirming that the construction extends naturally from the scalar to the vector-valued setting.
The \ac{FEM} results serve as a reference and show the expected behavior, that is, all three coarse space variants produce stable condition numbers without any moment correction, since all \ac{FEM} \ac{DOFs} are nodal.

\begin{table}[tbp]
\centering
\caption{Condition numbers and PCG iteration counts for the linear elasticity problem on a CVT mesh with $256$ subdomains partitioned by METIS.}
\label{tab:elasticity_degree_cvt}
\begin{adjustbox}{width=\textwidth}
\begin{tabular}{|c|l||r||r|r||r|r||r|r||r|r||r|r|}
\cline{3-13}
\multicolumn{2}{c|}{} & Coarse & \multicolumn{2}{c||}{$k=1$} & \multicolumn{2}{c||}{$k=2$} & \multicolumn{2}{c||}{$k=3$} & \multicolumn{2}{c||}{$k=4$} & \multicolumn{2}{c|}{$k=5$} \\
\cline{3-13}
\multicolumn{2}{c|}{} & & $\kappa$ & it & $\kappa$ & it & $\kappa$ & it & $\kappa$ & it & $\kappa$ & it \\
\cline{3-13}
\noalign{\vskip 0.6ex}
\cline{1-13}
\multirow{7}{*}{\STAB{\rotatebox[origin=c]{90}{VEM}}}
& ASM (1-lv)           & \multicolumn{1}{c||}{--} & \multicolumn{2}{c||}{N.C.} & \multicolumn{2}{c||}{N.C.} & \multicolumn{2}{c||}{N.C.} & \multicolumn{2}{c||}{N.C.} & \multicolumn{2}{c|}{N.C.} \\
& GDSW$^{\rm std}$     & 15560 & $16.76$ & 35 & $22.83$ & 41 & \multicolumn{2}{c||}{N.C.} & \multicolumn{2}{c||}{N.C.} & \multicolumn{2}{c|}{N.C.} \\
& GDSW$^{*,\rm std}$   & 9552  & $16.91$ & 35 & $22.47$ & 41 & \multicolumn{2}{c||}{N.C.} & \multicolumn{2}{c||}{N.C.} & \multicolumn{2}{c|}{N.C.} \\
& RGDSW$^{\rm std}$    & 3306  & $18.46$ & 38 & $25.34$ & 44 & \multicolumn{2}{c||}{N.C.} & \multicolumn{2}{c||}{N.C.} & \multicolumn{2}{c|}{N.C.} \\
& GDSW                 & 15560 & $16.76$ & 35 & $22.83$ & 41 & $26.79$             &   45 & $29.44$            &   47 & $32.1$ &   49 \\
& GDSW$^*$             & 9552  & $16.91$ & 35 & $22.47$ & 41 & $26.42$             &   45 & $28.95$            &   47 & $31.7$ &   49 \\
& RGDSW                & 3306  & $18.46$ & 38 & $25.34$ & 44 & $31.71$             &   51 & $35.47$            &   54 & $39.3$ &   58 \\
\hline
\end{tabular}
\end{adjustbox}
\end{table}

\Cref{tab:elasticity_degree_cvt} reports the same comparison for the CVT polyhedral mesh.
The results confirm that the moment-based construction is also effective on general polyhedral meshes with irregular partitions: the standard variants again fail for $k\ge 3$, while the corrected variants maintain bounded condition numbers for all tested degrees, in agreement with the CUBE mesh results.

\section{Conclusion}
\label{sec:conclusions}

This work extends GDSW-type overlapping Schwarz preconditioning to 3D high-order \ac{FEM} and \ac{VEM} discretizations of the Poisson equation, the scalar diffusion equation, and the equations of linear elasticity. In particular, we highlight the role of face moments in the coarse space construction, and we show that the same approach carries over componentwise to the vector-valued elasticity setting.
The numerical tests show that this modification is necessary for $k\ge3$ and remains effective on both structured and polyhedral meshes, for both model problems.

For Poisson's problem, for the tested polynomial degrees and decompositions, all GDSW variants produce stable iteration behavior, while RGDSW gives the smallest coarse space. The parallel timings demonstrate that the proposed implementation scales well to large core counts, while the RGDSW coarse problem remains cheap to solve without a corresponding loss in iteration robustness.

For the heterogeneous diffusion configuration depicted in \Cref{fig:channels}, GDSW and GDSW$^*$ remain robust while RGDSW degrades; for general coefficient distributions, robustness can be addressed using adaptive strategies.

\section*{Acknowledgments}
This project has received funding from the Deutsche Forschungsgemeinschaft (DFG) as part of the Forschungsgruppe (Research Unit) 5134 ``Solidification Cracks During Laser Beam Welding -- High Performance Computing for High Performance Processing'' under DFG project number 434946896. The authors gratefully acknowledge the scientific support and HPC resources provided by the Erlangen National High Performance Computing Center (NHR@FAU) of the Friedrich-Alexander-Universit\"at Erlangen-N\"urnberg (FAU) under the NHR project k109be10. NHR funding is provided by federal and Bavarian state authorities. NHR@FAU hardware is partially funded by the German Research Foundation (DFG) - 440719683. We would like to thank Franco Dassi from the University of Milano-Bicocca and the other developers of the \texttt{VEM++} library.

\section*{Author contributions and AI usage} All authors are listed in alphabetical order. The contributions are specified below using the CRediT taxonomy. TB: Conceptualization, Methodology, Data Curation, Software, Writing - original draft \& editing. AK: Conceptualization, Methodology, Writing - review \& editing. ML: Conceptualization, Methodology, Writing - review \& editing. AW: Conceptualization, Methodology, Writing - original draft \& editing.
AI models were used to assist with formatting and improving the clarity of portions of the text, as well as to assist with the development and debugging of the \texttt{FE2TI}-\texttt{VEM++} interface, specifically the analysis of DOF orderings and mappings within the \texttt{VEM++} library. The authors assume full responsibility for all content.

\bibliographystyle{siamplain}
\bibliography{refer}

\end{document}